\documentclass{amsart}

\usepackage{amscd,amsthm,amssymb,xspace} %
\usepackage{enumerate}        %
\numberwithin{equation}{section}

\newenvironment{roenumerate}{\begin{enumerate}[\ulp\itshape i\urp]}{\end{enumerate}}

\newtheorem{thm}[equation]{Theorem} %
\newtheorem{theorem}[equation]{Theorem}
\newtheorem{prop}[equation]{Proposition}
\newtheorem{proposition}[equation]{Proposition}
\newtheorem{cor}[equation]{Corollary}
\newtheorem{corollary}[equation]{Corollary}
\newtheorem{lemma}[equation]{Lemma}

\theoremstyle{definition}  %
\newtheorem{defn}[equation]{Definition}
\newtheorem{definition}[equation]{Definition}
\newtheorem{example}[equation]{Example}

\DeclareMathOperator{\im}{im}
\DeclareMathOperator{\Hom}{Hom}
\DeclareMathOperator{\Ext}{Ext}
\DeclareMathOperator{\Tor}{Tor}
\DeclareMathOperator{\colim}{colim}

\DeclareMathOperator{\cofin}{cofin}

\newcommand{\llp}{left lifting property with respect to\xspace}
\newcommand{\rlp}{right lifting property with respect to\xspace}

\newcommand{\Z}{\mathbb{Z}}
\newcommand{\Q}{\mathbb{Q}}

\newcommand{\Mod}{\text{-mod}}

\newcommand{\RMod}{\ensuremath{R}\Mod\xspace}
\newcommand{\SMod}{\ensuremath{S}\Mod\xspace}

\newcommand{\boxprod}{\mathbin{\square}}

\newcommand{\ulp}{\textup{(}}
\newcommand{\urp}{\textup{)}}

\newcommand{\mathcolon}{\colon\,}
\newcommand{\mc}{\mathcolon}
\newcommand{\ra}{\xrightarrow{}}
\newcommand{\lra}{\xrightarrow{}}
\newcommand{\llra}{\xrightarrow}
\newcommand{\tensor}{\otimes}
\newcommand{\Ho}{\text{Ho}}

\newcommand{\Ch}[1]{\text{Ch} (#1)}
\newcommand{\cat}[1]{\mathcal{#1}}
\newcommand{\cA}              {{\cat A}}
\newcommand{\sA}              {{s\cA}}
\newcommand{\cB}              {{\cat B}}
\newcommand{\cC}              {{\cat C}}
\newcommand{\cD}              {{\cat D}}
\newcommand{\cE}              {{\cat E}}
\newcommand{\cF}              {{\cat F}}
\newcommand{\cG}              {{\cat G}}
\newcommand{\cH}              {{\cat H}}

\newcommand{\cK}                {{\cat K}}
\newcommand{\cM}              {{\cat M}}
\newcommand{\cP}              {{\cat P}}
\newcommand{\cS}              {{\cat S}}

\newcommand{\Ab}                {{\cat Ab}}

\newcommand{\A}                 {{\cat A}}
\newcommand{\ChA}               {\Ch{\cA}}
\newcommand{\ChR}               {\Ch{R}}
\newcommand{\K}                 {\cK}
\newcommand{\KA}                {\K(\cA)}
\newcommand{\KR}                {\K(R)}
\newcommand{\D}                 {\cD}
\newcommand{\DA}                {\D(\cA)}
\newcommand{\DR}                {\D(R)}
\newcommand{\DC}                {\D_{\cC}}
\newcommand{\DP}                {\D_{\cP}}

\newcommand{\dfn}{\textbf} %
\newcommand{\mdfn}[1]{\dfn{\mathversion{bold}#1}} %

\DeclareMathOperator{\Ex}     {Ex}

\newcommand{\op}               {\text{op}}

\newcommand{\Proj}              {\cP}

\newcommand{\Epi}               {\cE}

\DeclareMathOperator{\Pure}     {P}
\newcommand{\PExt}              {\Pure\!\Ext}

\newcommand{\Iinj}              {\textup{$I$-inj}}
\newcommand{\Jinj}              {\textup{$J$-inj}}
\newcommand{\Icof}              {\textup{$I$-cof}}
\newcommand{\Jcof}              {\textup{$J$-cof}}
\newcommand{\Icell}              {\textup{$I$-cell}}
\newcommand{\Jcell}              {\textup{$J$-cell}}

\newcommand{\st}{\,\mid\,}
\newcommand{\ie}{\emph{i.e.,}\xspace}

\newcommand{\comment}[1]{}

\begin{document}
 
\title{Quillen model structures for relative homological algebra}

\date{February 24, 2013}

\author{J. Daniel Christensen}
\thanks{The first author was supported in part by NSF grant DMS 97-29992.}
\address{Department of Mathematics \\ University of Western Ontario \\
London, ON N6A 5B7 \\ Canada}
\email{jdc@uwo.ca}

\author{Mark Hovey}
\thanks{The second author was supported in part by NSF grant DMS 99-70978.}
\address{Department of Mathematics \\ Wesleyan University
\\ Middletown, CT 06459 \\ USA}
\email{hovey@member.ams.org}

\thanks{Published as \emph{Math.\ Proc.\ Cambridge Philos.\ Soc.}\ \textbf{133} Part 2 (2002), 
261-293.
See last page for corrections made since publication.}



\begin{abstract}
An important example of a model category is the category of unbounded
chain complexes of $R$-modules,
which has as its homotopy category the derived category of the ring $R$.
This example shows that traditional homological algebra is encompassed
by Quillen's homotopical algebra.
The goal of this paper is to show that more general forms of homological 
algebra also fit into Quillen's framework.
Specifically, a projective class on a complete and cocomplete abelian category
$\cA$ is exactly the information needed to do homological
algebra in $\cA$.  
The main result is that, under weak hypotheses,
the category of chain complexes of objects of $\cA$ has a 
model category structure that reflects the homological algebra
of the projective class in the sense that it encodes the
$\Ext$ groups and more general derived functors.
Examples include the ``pure derived category'' of a ring $R$,
and derived categories capturing relative situations, including
the projective class for Hochschild homology and cohomology.
We characterize the model structures that are cofibrantly generated,
and show that this fails for many interesting examples.
Finally, we explain how the category of simplicial objects in 
a possibly non-abelian category can be equipped with a model category
structure reflecting a given projective class, and give examples
that include equivariant homotopy theory and bounded below 
derived categories.
\end{abstract}

\vspace*{-2.5ex}
\maketitle

\vspace{-1.5\baselineskip}
\enlargethispage{\baselineskip}
\tableofcontents
\vspace*{-1.9\baselineskip} %
\pagebreak

\comment{This version includes comments to the author!}

\comment{I updated the abstract and the intro, but they could probably
still use polishing.}

\section*{Introduction}\label{se:intro}

An important example of a model category is the category $\ChR$ of unbounded
chain complexes of $R$-modules,
which has as its homotopy category the derived category $\DR$ of the 
associative ring $R$.
The formation of a projective resolution is an example of cofibrant
replacement, traditional derived functors are examples of
derived functors in the model category sense, and $\Ext$ groups appear
as groups of maps in the derived category.
This example shows that traditional homological algebra is encompassed
by Quillen's homotopical algebra, and indeed this unification was one
of the main points of Quillen's influential work~\cite{qu:ha}.

The goal of this paper is to illustrate that more general forms of 
homological algebra also fit into Quillen's framework.
In any abelian category $\A$ there is a natural notion of ``projective
object'' and ``epimorphism.''  
However, it is sometimes useful to impose different definitions
of these terms.
If this is done in a way that satisfies some natural axioms,
what is obtained is a ``projective class,'' which is exactly 
the information needed to do homological algebra in $\A$.  
Our main result shows that for a wide variety of projective classes
(including all those that arise in examples)
the category of unbounded chain complexes of objects of $\A$ has a 
model category structure that reflects the homological algebra
of the projective class in the same way that ordinary homological
algebra is captured by the usual model structure on $\ChR$.

When $\A$ has enough projectives, the projective objects and
epimorphisms form a projective class.  Therefore the results
of this paper apply to traditional homological algebra as well.
Even in this special case, it is not a trivial fact that the
category of \emph{unbounded} chain complexes can be given a model 
category structure, and indeed Quillen restricted himself to the
bounded below case.
We know of three other written proofs that the category of unbounded
chain complexes is a model category~\cite{gr:mc,hi:haha,hovey-model}, which do the
case of $R$-modules, but this was probably known to others as well.

An important corollary of the fact that a derived category $\DA$ is the
homotopy category of a model category is that the group $\DA(X,Y)$ of maps
is a set (as
opposed to a proper class) for any two chain complexes $X$ and $Y$.
This is not the case in general, and much work on derived categories 
ignores this possibility.  
The importance of this point is that if one uses the morphisms in the
derived category to index constructions in other categories or to
define cohomology groups, one needs to know that the indexing class
is actually a set.
Recently, $\DA(X,Y)$ has been shown to be a set
under various assumptions on $\A$.
(See Weibel~\cite{we:iha} Remark 10.4.5, which credits Gabber, and
Exercise 10.4.5, which credits Lewis, May and Steinberger~\cite{lemast:esht}.
See also Kriz and May~\cite[Part III]{krma:oamm}.)
The assumptions that appear in the present paper are different from those that
have appeared before and the proof is somewhat easier (because
of our use of the theory of cofibrantly generated model categories),
so this paper may be of some interest even in this special case.

Another consequence of the fact that $\ChR$ is a model category is
the existence of resolutions coming from cofibrant and fibrant
approximations, and the related derived functors.  
Some of these are discussed in \cite{avfoha} and \cite{sp:ruc}.
We do not discuss these topics here, but just mention that
these resolutions are immediate once you have the model structure,
so our approach gives these results with very little work.

While our results include new examples of traditional homological
algebra, our focus is on more general projective classes.
For example, let $A$ be an algebra over a commutative ring $k$.
We call a map of $A$-bimodules a relative epimorphism if it
is split epic as a map of $k$-modules, and we call an $A$-bimodule
a relative projective if maps from it lift over relative epimorphisms.
These definitions give a projective class, and Theorem~\ref{thm-rings}
tells us that there is a model category, and therefore a derived category,
that captures the homological algebra of this situation.  For example,
Hochschild cohomology groups appear as $\Hom $ sets in this derived 
category (see Example~\ref{ex:bimodules}).

We also discuss pure homological algebra and construct the
``pure derived category'' of a ring.
Pure homological algebra has applications to phantom maps in
the stable homotopy category~\cite{chst:pmht} and in the (usual) derived
category of a ring~\cite{ch:itcpgs}, connections to
Kasparov KK-theory~\cite{sc:fineI}, and is actively studied
by algebraists and representation theorists.

In the last section we describe a model category structure
on the category of non-negatively graded chain complexes
that works for an arbitrary projective class on an abelian category, 
without any hypotheses.
More generally, we show that under appropriate hypotheses a projective
class on a possibly non-abelian category $\cA$ determines a model
category structure on the category of simplicial objects in $\cA$.
As an example, we deduce that the category of equivariant simplicial
sets has various model category structures.

We now briefly outline the paper.  In Section~\ref{se:pc} we
give the axioms for a projective class and mention many examples
that will be discussed further in Subsections~\ref{subse:adjoint-examples} 
and~\ref{subse:pure}.
In Section~\ref{se:model-general} we describe the desired model
structure coming from a projective class and state our main theorem, 
which says that the model structure exists as long as cofibrant 
replacements exist.  We also give two hypotheses that each
imply the existence of cofibrant replacements.  
The first hypothesis handles situations coming from adjoint pairs,
and is proved to be sufficient in Section~\ref{se:adjoints}, where
we also give many examples involving relative situations.
The second hypothesis deals with projective classes that 
have enough small projectives and is proved
to be sufficient in Section~\ref{se:model-small-projs}.
In Section~\ref{se:model-set-small-projs} we prove that the model structure
that one gets is cofibrantly generated if and only if there
is a \emph{set} of enough small projectives.  We do this
using the recognition lemma for cofibrantly generated categories, 
which is recalled in Subsection~\ref{subse:cof}.
This case is proved from scratch, independent of the main result in 
Section~\ref{se:model-general}, since the proof is not long.
In Subsection~\ref{subse:pure} we give two examples, the traditional
derived category of $R$-modules and the pure derived category.
We describe how the two relate and why the pure derived category is 
interesting.
In the final section we discuss the bounded below case,
which works for \emph{any} projective class, and describe a result for
simplicial objects in a possibly non-abelian category.

We thank Haynes Miller for asking the question that led to this paper
and Haynes Miller and John Palmieri for fruitful and
enjoyable discussions.

\setcounter{section}{-1}
\section{Notation and conventions}\label{se:notation}

We make a few blanket assumptions.  With the exception of 
Section~\ref{se:simp}, $\cA$ and $\cB$ will denote abelian categories.  
We will assume that our abelian categories are bicomplete; this
assumption is stronger than strictly necessary, but it simplifies the
statements of our results.
For any category $\cC$, we write $\cC(A,B)$ for the set of maps
from $A$ to $B$ in $\cC$.

We write $\ChA$ for the category of unbounded chain complexes of
objects of $\cA$ and degree zero chain maps.
To fix notation, assume that the differentials lower degree.
For an object $X$ of $\ChA$, define $Z_n X := \ker(d \mc X_n \ra X_{n-1})$ and
$B_n X := \im(d \mc X_{n+1} \ra X_n)$, and write $H_n X$ for the quotient.
A map inducing an isomorphism in $H_{n}$ for all $n$ is a \dfn{quasi-isomorphism}.
The \dfn{suspension} $\Sigma X$ of $X$ has $(\Sigma X)_{n} = X_{n-1}$
and $d_{\Sigma X} = -d_{X}$.
The functor $\Sigma$ is defined on morphisms by $(\Sigma f)_n = f_{n-1}$.
Given a map $f \mc X \ra Y$ of chain complexes, the \mdfn{cofibre
of $f$} is the chain complex $C$ with $C_{n} = Y_{n} \oplus X_{n-1}$
and with differential $d(y,x) = (dy+fx,-dx)$.
There are natural maps $Y \ra C \ra \Sigma X$, and the
sequence $X \ra Y \ra C \ra \Sigma X$ is called a
\dfn{cofibre sequence}.

Two maps $f, g \mc X \ra Y$ are \dfn{chain homotopic} if there
is a collection of maps $s_{n} \mc X_{n} \ra Y_{n+1}$ such that
$f - g = d s + s d$.
We write $[X, Y]$ for the chain homotopy classes of maps from
$X$ to $Y$ and $\KA$ for the category of chain complexes and
chain homotopy classes of maps.
Two complexes are \dfn{chain homotopy equivalent} if they are
isomorphic in $\KA$, and a complex is \dfn{contractible} if
it is chain homotopy equivalent to $0$.
$\KA$ is a triangulated category, with triangles the sequences
chain homotopy equivalent to the cofibre sequences above.
The functors $H_{n}(-)$, $[X,-]$ and $[-,Y]$ are defined on
$\KA$ and send triangles to long exact sequences.

For $P$ in $\cA$, $D^{k}P$ denotes the (contractible) complex such that
$(D^{k}P)_{n}=P$ if $n=k$ or $n=k-1$ but $(D^{k}P)_{n}=0$ for other
values of $n$, and whose differential is the identity in degree $k$.
The functor $D^{k}$ is left adjoint to the functor $X \mapsto X_{k}$,
and right adjoint to the functor $X \mapsto X_{k-1}$.  The path complex
$PX$ of a complex $X$ is the contractible complex such that $(PX)_{n} =
X_{n} \oplus X_{n+1}$, where $d(x,y)=(dx, x-dy)$.

We assume knowledge of the basics of model categories, for 
which~\cite{dwsp:htmc} is an excellent reference.
We use the definition of model category that requires that
the category be complete and cocomplete, and that the
factorizations be functorial.
Correspondingly, when we say that cofibrant replacements
exist, we implicitly mean that they are functorial.

\section{Projective classes}\label{se:pc}

\subsection{Definition and some examples}

In this subsection we explain the notion of a projective class, which is
the information necessary in order to do homological algebra.
Intuitively, a projective class is a choice of which sort of 
``elements'' we wish to think about.
In this section we focus on the case of an abelian category, but
this definition works for any pointed category with kernels.

The elements of a set $X$ correspond bijectively to the
maps from a singleton to $X$, and
the elements of an abelian group $A$ correspond bijectively to
the maps from $\Z$ to $A$.
Motivated by this, we 
call a map $P \ra A$ in any category a \mdfn{$P$-element of $A$}.
If we don't wish to mention $P$, we call such a map a \mdfn{generalized
element of $A$}.
A map $A \ra B$ in any category is determined by what it does on
generalized elements.
If $\Proj$ is a collection of objects, then a \mdfn{$\Proj$-element} means
a $P$-element for some $P$ in $\Proj$.

Let $\A$ be an abelian category.  A map $B \ra C$
is said to be \mdfn{$P$-epic} if it induces a surjection of
$P$-elements, that is, if the induced map $\A(P,B) \lra \A(P,C)$
is a surjection of abelian groups.
The map $B \ra C$ is \mdfn{$\Proj$-epic} if it is $P$-epic for all
$P$ in $\Proj$.

\begin{defn}\label{de:pc}
A \dfn{projective class} on $\A$ is a collection $\Proj$ of objects
of $\A$ and a collection $\Epi$ of maps in $\A$ 
such that 
\begin{roenumerate}
\item $\Epi$ is precisely the collection of all $\Proj$-epic maps;
\item $\Proj$ is precisely the collection of all objects $P$ such that each
      map in $\Epi$ is $P$-epic;
\item for each object $B$ there is a map $P \ra B$ in $\Epi$ with $P$ in 
      $\Proj$.
\end{roenumerate}
When a collection $\Proj$ is part of a projective class $(\Proj,\Epi)$,
the projective class is unique, and so we say that $\Proj$ determines
a projective class or even that $\Proj$ is a projective class.
An object of $\Proj$ is called a \mdfn{$\Proj$-projective}, or, if
the context is clear, a \dfn{relative projective}.
\end{defn}

A sequence 
\[
A \lra B \lra C
\]
is said to be \mdfn{$P$-exact} if the composite $A \ra C$ is zero and
\[
\A(P,A) \lra \A(P,B) \lra \A(P,C)
\]
is an exact sequence of abelian groups.
The latter can be rephrased as the condition that $A \ra B \ra C$
induces an exact sequence of $P$-elements.
A \mdfn{$\Proj$-exact sequence} is one that is $P$-exact for all
$P$ in $\Proj$.

\begin{example}\label{ex:cp}
For an associative ring $R$,
let $\A$ be the category of left $R$-modules, let $\Proj$ be the collection
of all summands of free $R$-modules and let $\Epi$ be the collection of
all surjections of $R$-modules.  Then $\Epi$ is precisely the
collection of $\Proj$-epimorphisms, and $\Proj$ is a projective
class.
The $\Proj$-exact sequences are the usual exact sequences.
\end{example}

Example~\ref{ex:cp} is a \dfn{categorical} projective class in
the sense that the $\Proj$-epi\-mor\-phisms are just the epimorphisms
and the $\Proj$-projectives are the categorical projectives,
\ie those objects $P$ such that maps from $P$ lift through epimorphisms.

Here are two examples of non-categorical projective classes.

\begin{example}\label{ex:trivial-projective-class}
If $\A$ is any abelian category,
$\Proj$ is the collection of all objects, 
and $\Epi$ is the collection of all split epimorphisms $B \ra C$,
then $\Proj$ is
a projective class.  It is called the \dfn{trivial projective class}.
A sequence $A \ra B \ra C$ is $\Proj$-exact if and only if
$A \ra \ker(B \ra C)$ is split epic.
\end{example}

\begin{example}\label{ex:pp}
Let $\A$ be the category of left $R$-modules, as in Example~\ref{ex:cp}.
Let $\Proj$ consist of all summands of sums of finitely presented
modules and define $\Epi$ to consist of all $\Proj$-epimorphisms.
Then $\Proj$ is a projective class.
A sequence is $\Proj$-exact iff it is exact after
tensoring with every right module.
\end{example}

Examples~\ref{ex:cp} and~\ref{ex:pp} will be discussed further in
Subsection~\ref{subse:pure}.
Example~\ref{ex:trivial-projective-class} is important because
many interesting examples are ``pullbacks'' of this
projective class (see Subsection~\ref{subse:pullbacks}).

Let $\Proj $ be a projective class.
If $\cS$ is a subcollection of $\Proj$ (not necessarily a set),
and if a map is $\cS $-epic iff it is $\Proj$-epic, 
then we say that $\Proj$ is \mdfn{determined by $\cS $}
and that $\cS $ is a collection of \dfn{enough projectives}.
Some projective classes, such as Examples~\ref{ex:cp} and \ref{ex:pp}, 
are determined by a set, and the lemma below shows that any set of 
objects determines a projective class.
The trivial projective class is sometimes not determined by a set
(see Subsection~\ref{subse:not-det-by-a-set}).

\begin{lemma}\label{le:det-by-a-set}
Suppose $\cF$ is any \emph{set} of objects in an abelian category with
coproducts.
Let $\Epi$ be the collection of $\cF$-epimorphisms and
let $\Proj$ be the collection of all objects $P$ such that every
map in $\Epi$ is $P$-epic.
Then $\Proj$ is the collection of retracts of coproducts of objects
of $\cF$ and $(\Proj,\Epi)$ is a projective class.
\end{lemma}

\begin{proof}
Given an object $X$, let $P$ be the coproduct $\coprod F$ indexed
by all maps $F \ra X$ and all objects $F$ in $\cF$.  
The natural map $P \ra X$ is clearly an $\cF$-epimorphism.
Moreover, if $X$ is in $\Proj$, then this map is split epic,
and so $X$ is a retract of a coproduct of objects of $\cF$.
These two facts show that $(\Proj,\Epi)$ is a projective class.
\end{proof}

\subsection{Homological algebra}\label{subse:hom-alg}

A projective class is precisely the information needed to form
projective resolutions and define derived functors.
All of the usual definitions and theorems go through.
A \mdfn{$\Proj$-resolution} of an object $M$ is a $\Proj$-exact
sequence
\[
\cdots \lra P_{2} \lra P_{1} \lra P_{0} \lra M \lra 0
\]
such that each $P_{i}$ is in $\Proj$.
If $\cB$ is an abelian category and $T \mc \cA \ra \cB$ is an
additive functor, then the $n$th left derived functor of $T$
with respect to $\Proj$ is defined by 
$L_{n}^{\Proj}T(M) = H_{n}(T(P_{*}))$ where $P_{*}$ is a
$\Proj$-resolution of $M$.
One has the usual uniqueness of resolutions up to chain homotopy
and so this is well-defined.
 From a $\Proj$-exact sequence $0 \ra L \ra M \ra N \ra 0$ one
gets a long exact sequence involving the derived functors.
The abelian groups $\Ext_{\Proj}^{n}(M,N)$ can be defined
in the usual two ways, as equivalence classes of $\Proj$-exact
sequences $0 \ra N \ra L_{1} \ra \cdots \ra L_{n} \ra M \ra 0$,
or as $L_{n}^{\Proj}T(M)$ where $T(-) = \cA(-,N)$.

For further details and useful results we refer the reader to the
classic reference~\cite{eimo:frha}.

\subsection{Pullbacks}\label{subse:pullbacks}

A common setup in relative homological algebra is the following.  
We assume we have a functor $U \mc \cA \ra \cB$ of abelian categories, 
together with a left adjoint $F \mc \cB \ra \cA$.  
Then $U$ and $F$ are additive, $U$ is left exact and $F$ is right exact.  

If $(\Proj', \Epi')$ is a projective class on $\cB$, we define
$\Proj := \{\text{retracts of $FP$ for $P$ in $\Proj'$}\}$ and
$\Epi := \{\text{$B \ra C$ such that $UB \ra UC$ is in $\Epi'$}\}$.
Then one can easily show that $(\Proj, \Epi)$ is a 
projective class on $\cA$ and that 
a sequence is $\Proj$-exact if and only if it is sent to a 
$\Proj'$-exact sequence by $U$.
$(\Proj, \Epi)$ is called the \dfn{pullback}
of $(\Proj', \Epi')$ along the right adjoint $U$.

The most common case is when $(\Proj', \Epi')$ is the trivial
projective class (see Example~\ref{ex:trivial-projective-class}).
Then for any $M$ in $\cA$ the counit $F U M \ra M$ is
a $\Proj$-epimorphism from a $\Proj$-projective.

\begin{example}\label{ex:change-of-rings}
Let $R \ra S$ be a map of associative rings.  Write \RMod and \SMod for the
categories of left $R$- and $S$-modules.
Consider the forgetful functor $U \mc \SMod \ra \RMod$  and its 
left adjoint $F$ that sends an $R$-module $M$ to $S \tensor_{R} M$.
The pullback along $U$ of the trivial projective class on $\RMod$
gives a projective class $\Proj$.
The $\Proj$-projectives are the $S$-modules $P$ such that the natural
map $S \tensor_{R} P \ra P$ is split epic as a map of $S$-modules.
The $\Proj$-epimorphisms are the $S$-module maps that are split
epic as maps of $R$-modules.
\end{example}

\begin{example}\label{ex:change-of-rings-inj}
As above, let $R \ra S$ be a map of rings.
The forgetful functor $U \mc \SMod \ra \RMod$ has a
right adjoint $G$ that sends an $R$-module $M$ to the $S$-module $\RMod(S,M)$.
We can pullback the trivial injective class along $U$ to get an 
injective class on $\SMod$.
(An injective class is just a projective class on the opposite category.)
The relative injectives are the $S$-modules $I$ such that the natural
map $I \ra \RMod(S,I)$ is split monic as a map of $S$-modules,
and the relative monomorphisms are the $S$-module maps that
are split monic as maps of $R$-modules.
\end{example}

We investigate these examples in detail in Section~\ref{se:adjoints}.

\subsection{Strong projective classes}\label{subse:strong}

In Section~\ref{se:adjoints} we will focus on projective classes
that are the pullback of a trivial projective class along a 
right adjoint.
In this subsection we describe a special property that these projective
classes have.

\begin{defn}
A projective class $\Proj$ is \dfn{strong} if for each $\Proj$-projective $P$
and each $\Proj$-epimorphism $M \ra N$, the surjection
$\cA(P,M) \ra \cA(P,N)$ of abelian groups is \emph{split} epic.
\end{defn}

It is clear that a trivial projective class is strong, and that
the pullback of a strong projective class is strong.  
The importance of strong projective classes comes from the
following lemma.

\begin{lemma}\label{le:strong}
The following are equivalent for a projective class $\Proj$ on an
abelian category $\cA$:
\begin{roenumerate}
\item $\Proj$ is strong.
\item For each complex $C$ in $\ChA$, if the complex $\cA(P,C)$ in $\Ch{\Ab}$
      has trivial homology for each $P$ in $\Proj$,
      then it is contractible for each $P$ in $\Proj$.
\item For each map $f$ in $\ChA$, if the map $\cA(P,f)$ in $\Ch{\Ab}$ 
      is a quasi-isomorphism for each $P$ in $\Proj$, 
      then it is a chain homotopy equivalence for each $P$ in $\Proj$.
\end{roenumerate}
\end{lemma}

Here $\cA (P,C)$ denotes the chain complex with $\cA (P,C_{k})$
in degree $k$.
We also use the notation from Section~\ref{se:notation}.

\begin{proof}
(i) $\implies$ (ii):  Assume $\Proj$ is a strong projective class and 
let $C$ be a complex in $\ChA$ such that $\cA(P,C)$ has trivial
homology for each $\Proj$-projective $P$.
Then for each $k$ and each $P$ we have a short exact sequence
\[
0 \lra \cA(P, Z_{k}C) \lra \cA(P, C_{k}) \lra \cA(P, Z_{k-1}C) \lra 0
\]
of abelian groups.
Because the projective class is strong, the sequence is split.
This implies that the complex $\cA(P,C)$ is isomorphic to
$\oplus_{k} D^{k+1} \cA(P, Z_{k}C)$, and in particular that
it is contractible.

(ii) $\implies$ (iii):  Let $f \mc X \ra Y$ be a map in $\ChA$ such that
$\cA(P,X) \ra \cA(P,Y)$ is a quasi-isomorphism for each $P$ in $\Proj$,
and let $C$ be the cofibre of $f$.
Then by the long exact sequence, $\cA(P,C)$ has trivial homology for
each $P$ in $\Proj$.  By (ii) this complex is contractible.
This implies that $\cA(P,X) \ra \cA(P,Y)$ is a chain homotopy equivalence.

(iii) $\implies$ (i): Let $M \ra N$ be a $\Proj$-epimorphism with kernel $L$.
Then the complex $L \ra M \ra N$ has trivial homology after applying 
$\cA(P,-)$, for each $\Proj$-projective $P$.  
By (iii), it is contractible after applying $\cA(P,-)$, for each $P$.  
In particular, $\cA(P,M) \ra \cA(P,N)$ is split epic.
\end{proof}

\section{The relative model structure}\label{se:model-general}

The object of this section is to construct a Quillen model
structure on the category $\ChA$ of chain complexes over
$\A$ that reflects a given projective class $\Proj$ on $\A$.

If $X$ is a chain complex, we write
$\A(P,X)$ for the chain complex that has the abelian group
$\A(P,X_{n})$ in degree $n$.  This is the chain complex of
$P$-elements of $X$.  

\begin{definition}\label{de:P-structure}
A map $f\mc X \lra Y$ in $\ChA$ is a \mdfn{$\Proj$-equivalence} 
if $\A(P,f)$ is a quasi-isomorphism in $\Ch{\Z}$ for each $P$ in $\Proj$.
The map $f$ is a \mdfn{$\Proj$-fibration} 
if $\A(P,f)$ is a surjection for each $P$ in $\Proj$.
The map $f$ is a \mdfn{$\Proj$-cofibration} 
if $f$ has the \llp all maps that are both $\Proj$-fibrations 
and $\Proj$-equivalences \ulp the \mdfn{$\Proj$-trivial fibrations}\urp.
\end{definition}

The motivation for this definition is that it implies that a
complex
\[
\cdots \lra P_{2} \lra P_{1} \lra P_{0} \lra 0 \lra \cdots
\]
equipped with an augmentation $P_{0} \ra M$ to an object $M$
is a cofibrant replacement if and only if it is a $\Proj$-resolution
in the sense of Subsection~\ref{subse:hom-alg}.
This implies that if $M$ and $N$ are objects of $\cA$ thought
of as complexes concentrated in degree zero, then
$\Ext_{\Proj}^{n}(M,N)$ can be identified with maps
from $\Sigma^{n} M$ to $N$ in the homotopy category of the
model category $\Ch{\cA}$.  
This will be described in more detail in Subsection~\ref{subse:properties}.

The main goal of this section is then to prove the following theorem. 

\begin{theorem}\label{thm-rings}
Suppose $\Proj$ is a projective class on the abelian category $\A$.
Then the category $\Ch{\cA}$, together with the $\Proj$-equivalences, the
$\Proj$-fibrations, and the $\Proj$-cofibrations, forms a Quillen
model category if and only if cofibrant replacements exist.
When the model structure exists, it is proper.
Cofibrant replacements exist in each of the following cases:
\begin{itemize}
\item [A:] $\Proj$ is the pullback of the trivial projective class
           along a right adjoint that preserves countable sums.
\item [B:] There are enough $\kappa$-small $\Proj$-projectives for
           some cardinal $\kappa$, and $\Proj$-resolutions can
           be chosen functorially.
\end{itemize}
\end{theorem}
The words ``enough'' and ``$\kappa$-small'' will be explained in
Section~\ref{se:model-small-projs}.
\comment{Mention that categorical projective class on opposite of
a Grothendieck category is an example too?  (If so, ref. Tibor.)
Are there any other known examples that don't fit into one of 
these cases?}

We call this structure the \mdfn{relative model structure}.  We point
out that we are using the modern definition of model
category~\cite{dwhika, hovey-model}, so our factorizations will be
functorial.  
Correspondingly, we require our cofibrant replacements to be functorial 
as well.
This theorem requires our blanket assumption that abelian
categories are bicomplete.

In Subsection~\ref{subse:properties} we will describe further
properties of these model structures, including conditions
under which they are monoidal.
In Section~\ref{se:model-set-small-projs} we show that if there
is a \emph{set} of enough small projectives, then the model structure
is cofibrantly generated.
On the other hand, we show in Subsection~\ref{subse:not-det-by-a-set}
that model categories coming from Case A are generally not
cofibrantly generated.

In Subsection~\ref{subse:why} we explain why it is no loss of
generality to start with a projective class $(\Proj,\Epi)$, rather than 
just an arbitrary class $\Proj$ of test objects.

\begin{proof}
Some of the properties necessary for a model category are evident from
the definitions.  It is clear that $\Ch{\cA}$ is bicomplete, since
$\cA$ is so.  Also, $\Proj$-equivalences have the two out
of three property, and $\Proj$-equivalences,
$\Proj$-fibrations, and $\Proj$-cofibrations are closed under
retracts.  Furthermore, $\Proj$-cofibrations have the \llp
$\Proj$-trivial fibrations, by definition.  It remains to show that
$\Proj$-trivial cofibrations have the \llp $\Proj$-fibrations, and
that the two factorization axioms hold.  The remaining lifting property
will be proved in Proposition~\ref{prop-trivial-cofib}, and the two
factorization axioms will be proved in
Propositions~\ref{prop-cofib-triv-fib} and~\ref{prop-triv-cofib-fib},
assuming that cofibrant replacements exist.
Properness will be proved in Proposition~\ref{prop-proper}, which
also defines the term.

That cofibrant replacements exist in cases A and B will be proved in 
Sections~\ref{se:adjoints} and 
\ref{se:model-small-projs}, respectively.
\end{proof}

Note that our Theorem~\ref{thm-rings} will also apply when we
have an injective class, that is, a projective class on $\A^{\op}$,
by dualizing the definition of the model structure.

We continue the proof of Theorem~\ref{thm-rings} with
a lemma that gives us a simple test of the lifting property.  
We use the notation from Section~\ref{se:notation}.

\begin{lemma}\label{lem-test-lift}
Suppose $p\mathcolon X\xrightarrow{}Y$ is a $\Proj$-fibration with
kernel $K$, and $i\mathcolon A\xrightarrow{}B$ is a degreewise split
inclusion whose cokernel $C$ is a complex of relative projectives.  If
every map $C\xrightarrow{}\Sigma K$ is chain homotopic to $0$, then $i$
has the \llp $p$.
\end{lemma}

\begin{proof}
We can write $B_{n}\cong A_{n}\oplus C_{n}$, where the differential is
defined by $d(a,c)=(da+\tau c,dc)$ (we use the element notation for
convenience, but it is not strictly necessary), and $\tau \mathcolon
C_{n}\xrightarrow{}A_{n-1}$ can be any family of maps such that $d\tau
+\tau d=0$.  Suppose we have a commutative square as below.  
\[
\begin{CD}
A @>f>> X \\
@ViVV @VVpV \\
B @>>g> Y
\end{CD}
\]
In terms of the splitting $B_{n}\cong A_{n}\oplus C_{n}$, we have
$g(a,c)=pf(a)+\alpha (c)$, where the family $\alpha_{n} \mathcolon
C_{n}\xrightarrow{}Y_{n}$ satisfies $d\alpha =pf\tau +\alpha d$.  We are
looking for a map $h\mathcolon B\xrightarrow{}X$ making the diagram
above commute.  In terms of the splitting, this means we are looking for
a family of maps $\beta_{n} \mathcolon C_{n}\xrightarrow{}X_{n}$ such
that $p\beta =\alpha $ and $d\beta =f\tau +\beta d$.  Since $C_{n}$ is
relatively projective, and $p$ is $\Proj$-epic, there is a map
$\gamma \mathcolon C_{n}\xrightarrow{}X_{n}$ such that $p\gamma =\alpha$.  
The difference $\delta =d\gamma -f\tau -\gamma d$ may not be zero,
but at least $p\delta =0$.  Let $j\mathcolon K\xrightarrow{}X$ denote
the kernel of $p$.  Then there is a map $F\mathcolon
C_{n}\xrightarrow{}K_{n-1}$ such that $jF=\delta $.  Furthermore, one
can check that $Fd=-dF$, so that $F\mathcolon C\xrightarrow{}\Sigma K$
is a chain map.  By hypothesis, $F$ is chain homotopic to $0$ by a map
$D\mathcolon C_{n}\xrightarrow{}K_{n}$, so that $Dd-dD=F$.  Define
$\beta =\gamma + jD$.  Then $\beta $ defines the desired lift, so $i$
has the \llp $p$.  
\end{proof}

Now we study the $\Proj$-cofibrations.  A complex $C$ is called
\mdfn{$\Proj$-cofibrant} if the map $0\xrightarrow{}C$ is a
$\Proj$-cofibration.
A complex $K$ is called \mdfn{weakly $\Proj$-contractible} if
the map $K \lra 0$ is a $\Proj$-equivalence, or, equivalently, if
all maps from a complex $\Sigma^{k} P$ consisting
of a relative projective concentrated in one degree to $K$ are
chain homotopic to $0$.

\begin{lemma}\label{lem-cofibrant}
A complex $C$ is $\Proj$-cofibrant if and only if each $C_{n}$ is 
relatively projective and every map from $C$ to a weakly 
$\Proj$-contractible complex $K$ is chain homotopic to $0$.
\end{lemma}

\begin{proof}
Suppose first that $C$ is $\Proj$-cofibrant.
If $M\xrightarrow{}N$ is a $\Proj$-epimorphism, then the
map $D^{n+1}M\xrightarrow{}D^{n+1}N$ is a $\Proj$-fibration.  It is
also a $\Proj$-equivalence, since it is in fact a chain homotopy
equivalence.  Since $C$ is $\Proj$-cofibrant, the map
\[
\Ch{\cA}(C,D^{n+1}M)\xrightarrow{}\Ch{\cA}(C,D^{n+1}N)
\]
is surjective. But this map is isomorphic to the map
$\cA(C_{n},M)\xrightarrow{}\cA(C_{n},N)$, so $C_{n}$ is
relatively projective.

If $K$ is weakly $\Proj$-contractible, then the natural map $PK \ra K$ is a 
$\Proj$-trivial fibration.  Since $C$ is $\Proj$-cofibrant, any map 
$C \ra K$ factors through $PK$, which means that it is chain homotopic 
to $0$.

The converse follows immediately from Lemma~\ref{lem-test-lift},
since the kernel of a $\Proj$-trivial fibration is 
weakly $\Proj$-contractible.
\end{proof}

\begin{proposition}\label{prop-cofibrations}
A map $i\mathcolon A\xrightarrow{}B$ is a $\Proj$-cofibration if and only
if $i$ is a degreewise split monomorphism with $\Proj$-cofibrant cokernel.
\end{proposition}

\begin{proof}
Suppose first that $i$ is a $\Proj$-cofibration with cokernel $C$.
Since $\Proj$-cofibrations are closed under pushouts, %
it is clear that $C$ is $\Proj$-cofibrant.  The map
$D^{n+1}A_{n}\xrightarrow{}0$ is a $\Proj$-fibration and a
$\Proj$-equivalence.  Since $i$ is a $\Proj$-cofibration, the map
$A\xrightarrow{}D^{n+1}A_{n}$ that is the identity in degree $n$ extends
to a map $B\xrightarrow{}D^{n+1}A_{n}$.  In degree $n$, this map defines
a splitting of $i_{n}$.

Conversely, suppose that $i$ is a degreewise split monomorphism and
the cokernel $C$ of $i$ is $\Proj$-cofibrant.  We need to show that
$i$ has the \llp any $\Proj$-trivial fibration $p\mathcolon
X\xrightarrow{}Y$.  But this follows by combining Lemmas~\ref{lem-test-lift} 
and~\ref{lem-cofibrant}.
\end{proof}

The next lemma provides a source of $\Proj$-cofibrant objects,
including the $\Proj$-cellular complexes.

\begin{defn}\label{de:cellular}
Call a complex $C$ \mdfn{purely $\Proj$-cellular} if it
is a colimit of a colimit-preserving diagram
\[
0 = C^{0} \lra C^{1} \lra C^{2} \lra \cdots
\]
indexed by an ordinal $\gamma$, such that for each $\alpha < \gamma$ the map
$C^{\alpha} \ra C^{\alpha+1}$ is degreewise split monic with cokernel 
a complex of relative projectives with zero differential.
By ``colimit-preserving'' we mean that for each limit ordinal
$\lambda < \gamma$, the map 
$\colim_{\alpha < \lambda} C^{\alpha} \ra C^{\lambda}$ is an isomorphism.
We say $C$ is \mdfn{$\Proj$-cellular} if it is a retract of a purely 
$\Proj$-cellular complex.
\end{defn}

\begin{lemma}\label{lem-cofibs}
\begin{enumerate}
\item [(a)] If $D$ in $\Ch{\cA}$ is a complex of relative projectives
with zero differential, then $D$ is $\Proj$-cofibrant.
\item [(b)] If $D$ in $\Ch{\cA}$ is a bounded below complex of
relative projectives, then $D$ is $\Proj$-cofibrant.
\item [(c)] If $D$ in $\Ch{\cA}$ is $\Proj$-cellular,
then $D$ is $\Proj$-cofibrant.
\end{enumerate}
\end{lemma}

\begin{proof}
(a) follows immediately from Lemma~\ref{lem-cofibrant}.

(b) Let $D$ be a bounded below complex of relative projectives and
write $D^{\leq n}$ for the truncation of $D$ that agrees with $D$ in
degrees $\leq n$ and is $0$ elsewhere.  Then the map 
$D^{\leq n} \ra D^{\leq n+1}$ is degreewise split monic and has a
$\Proj$-cofibrant cokernel (by (a)), so is a $\Proj$-cofibration
by Proposition~\ref{prop-cofibrations}.  
Since $D^{\leq n}=0$ for $n<<0$ and $0$ is $\Proj$-cofibrant, 
each $D^{\leq n}$ is $\Proj$-cofibrant.  Therefore, so is their 
colimit $D$.

(c) The proof is just a transfinite version of the proof of (b),
combined with the fact that a retract of a cofibrant object is cofibrant.
\end{proof}

We can now prove that the other lifting axiom holds.  

\begin{proposition}\label{prop-trivial-cofib}
A map $i\mathcolon A\xrightarrow{}B$ has the \llp $\Proj$-fibrations
if and only if $i$ is a $\Proj$-trivial cofibration.
\end{proposition}

\begin{proof}
Suppose first that $i$ has the \llp $\Proj$-fibrations.  Then $i$ is
a $\Proj$-cofibration, by definition, and the cokernel
$0\xrightarrow{}C$ also has the \llp $\Proj$-fibrations.  In
particular, since the map $PC\xrightarrow{}C$ is a $\Proj$-fibration,
$C$ is contractible.  Hence $i$ is a chain homotopy equivalence, and in
particular a $\Proj$-trivial cofibration.

Conversely, suppose that $i$ is a $\Proj$-trivial cofibration with
cokernel $C$.  By Lemma~\ref{lem-test-lift}, in order to show that $i$
has the \llp $\Proj$-fibrations, it suffices to show that every map
from $C$ to any complex $K$ is chain homotopic to $0$.  This is
equivalent to showing that $C$ is contractible.  Since $i$ is
degreewise split monic, for each relative projective $P$ there
is a long exact sequence 
$\cdots \ra [\Sigma^{k} P, A] \ra [\Sigma^{k} P, B] \ra 
 [\Sigma^{k} P, C] \ra \cdots .$
Since $i$ is a $\Proj$-equivalence, $[\Sigma^{k} P, C] = 0$ for each
relative projective $P$ and each $k$,
and so $PC\xrightarrow{}C$
is a $\Proj$-trivial fibration.  Since $C$ is $\Proj$-cofibrant,
the identity map of $C$ factors through $PC$, and so $C$ is
contractible.  
\end{proof}

Note that the proof shows that a $\Proj$-trivial cofibration is in
fact a chain homotopy equivalence.

Now we proceed to prove the factorization axioms, under the
assumption that we have cofibrant replacements.

\begin{proposition}\label{prop-cofib-triv-fib}
If every object $A$ has a cofibrant replacement
$q_{A} \mc QA \ra A$, then
every map in $\Ch{\cA}$ can be factored into a
$\Proj$-cofibration followed by a $\Proj$-trivial fibration.  
\end{proposition}

\begin{proof}
Suppose $f \mc A \ra B$ is a map in $\Ch{\cA}$.
Let $C$ be the cofibre of $f$, so $C = B \oplus \Sigma A$ with
$d(b,a) = (db + fa, -da)$, and 
let $E$ be the fibre of the composite $g \mc QC \ra C \ra \Sigma A$,
so $E = A \oplus QC$ with $d(a,q) = (da - gq, dq)$
(the desuspension of the cofibre).
Consider the diagram
\[
\begin{CD}
A @>i>> E @>>> QC @>g>> \Sigma A \\
@| @. @Vq_{C}VV @| \\
A @>>f> B @>>> C @>>> \Sigma A 
\end{CD}
\]
whose rows are triangles in $\K(\cA)$.
There is a natural fill-in map $p \mc E \ra B$ defined by
$p(a,q) = f(a) + \pi_{B} q_{C} q$, where $\pi_{B} \mc C \ra B$ is
the projection.
The map $p$ makes the left-hand square commute in $\Ch{\cA}$ and
the middle square commute in $\K(\cA)$ 
(with the chain homotopy $s(a,q) = (0,a)$).
The map $i \mc A \ra E$ is a $\Proj$-cofibration since it is 
degreewise split and its cokernel $QC$ is $\Proj$-cofibrant.
Furthermore, since 
$QC \ra C$ is degreewise $\Proj$-epic and $\pi_{B} \mc C \ra B$ is
degreewise split epic, it follows that $p$ is degreewise 
$\Proj$-epic.
Applying the functor $[\Sigma^{k} P, -]$ gives two long exact
sequences, and from the five-lemma one sees that $[\Sigma^{k} P, p]$ is 
an isomorphism when $P$ is a relative projective.
Thus $f = p i$ is the required factorization.
\end{proof}

\begin{proposition}\label{prop-triv-cofib-fib}
If every object $A$ has a cofibrant replacement
$q_{A} \mc QA \ra A$, then
every map in $\Ch{\cA}$ can be factored into a $\Proj$-trivial
cofibration followed by a $\Proj$-fibration.  
\end{proposition}

\begin{proof}
It is well-known that we can factor any map $f \mc A \ra B$ in $\Ch{\cA}$ 
into a degreewise split monomorphism that is also a 
chain homotopy equivalence, followed by a degreewise split epimorphism.  
Since every degreewise split epimorphism is a $\Proj$-fibration, 
we may as well assume $f$ is a degreewise split monomorphism and a 
chain homotopy equivalence.

In this case, we apply Proposition~\ref{prop-cofib-triv-fib} to factor
$f=pi$, where $p$ is a $\Proj$-trivial fibration and $i$ is a
$\Proj$-cofibration.  Since $f$ is a chain homotopy equivalence, $i$
must be a $\Proj$-trivial cofibration, and so the proof is complete.  
\end{proof}

Note that the factorizations constructed in
Proposition~\ref{prop-cofib-triv-fib} and
Proposition~\ref{prop-triv-cofib-fib} are both functorial in the map
$f$, since we are implicitly assuming that cofibrant replacement is
functorial.  

The homotopy category of $\ChA$, formed by inverting the
$\Proj$-equivalences, is called the \dfn{derived category} of $\A$
(with respect to $\Proj$).  It is denoted $\DA$, and a fundamental
result in model category theory asserts that $\DA(X,Y)$ is a set
for each $X$ and $Y$.
In Exercise 10.4.5 of \cite{we:iha}, Weibel outlines an argument
that proves that $\DA(X,Y)$ is a set when there are enough
(categorical) projectives, $\Proj$ is the categorical projective
class, and $\A$ satisfies AB5.  A connection between Weibel's hypotheses
and Case B is that if $\A$ has enough projectives that are
small with respect to all filtered diagrams in $\A$, then AB5 holds.
The smallness condition needed for our theorem is weaker than
this.  (See Section~\ref{se:model-small-projs} for the precise
hypothesis.)

\subsection{Properties of the relative model structure}\label{subse:properties}

In this subsection, we investigate some of the properties of the relative
model structure.  We begin by showing that the model category notions
of homotopy, derived functor, suspension and cofibre sequence agree
with the usual notions.  Then we study properness and monoidal
structure.  We discuss cofibrant generation in
Section~\ref{se:model-set-small-projs}.

We assume throughout that $\Proj$ is a projective class
on an abelian category $\cA$ such that the relative model structure 
on $\ChA$ exists.

We first show that the notion of homotopy determined from the
model category structure 
corresponds to the usual notion of chain homotopy.

\begin{defn}\textup{(\cite{dwsp:htmc} or \cite{qu:ha}.)}
If $M$ is an object in a model category $\cC$, a \dfn{good cylinder object}
for $M$ is an object $M \times I$ and a factorization
$M \coprod M \llra{i} M \times I \llra{p} M$ of the codiagonal map,
with $i$ a cofibration and $p$ a weak equivalence.
(Despite the notation, $M \times I$ is \emph{not} in general a product 
of $M$ with an object $I$.)
A \dfn{left homotopy} between maps $f, g \mc M \ra N$ is a map
$H \mc M \times I \ra N$ such that the composite $H i$ is equal
to $f \coprod g \mc M \coprod M \ra N$, for some good cylinder
object $M \times I$.
\end{defn}

The notion of \dfn{good path object} $N^{I}$ for $N$ is dual to that of 
good cylinder object and leads to the notion of \dfn{right homotopy}.  
The following standard result can be found in \cite[Section 4]{dwsp:htmc}, 
for example.

\begin{lemma}\label{le:homotopy}
For $M$ cofibrant and $N$ fibrant, two maps $f, g \mc M \ra N$ 
are left homotopic if and only if they are right homotopic,
and both of these relations are equivalence relations and
respect composition.
Moreover, if $M \times I$ is a fixed good cylinder object
for $M$, then $f$ and $g$ are left homotopic if and only if
they are left homotopic using $M \times I$;  similarly for
a fixed good path object.            \qed
\end{lemma}

Because of the lemma, for $M$ cofibrant and $N$ fibrant we have
a well-defined relation of \dfn{homotopy} on maps $M \ra N$.
Quillen showed that the homotopy category of $\cC$, which is by
definition the category of fractions formed by inverting the
weak equivalences, is equivalent to the category consisting
of objects that are both fibrant and cofibrant with morphisms
being homotopy classes of morphisms.

Now we return to the study of the model category $\ChA$.

\begin{lemma}\label{le:rh}
Let $M$ and $N$ be objects of $\ChA$ with $M$ $\Proj$-cofibrant.
Two maps $M \ra N$ are homotopic if and only if they are chain homotopic.
\end{lemma}

\begin{proof}
We construct a factorization $M \oplus M \ra M \times I \ra M$ of the
codiagonal map $M \oplus M \ra M$ in the following way.
Let $M \times I$ be the chain complex that has 
$M_{n} \oplus M_{n-1} \oplus M_{n}$ in degree $n$.  
We describe the differential by saying that it sends
a generalized element $(m,\bar{m},m')$ in $(M \times I)_{n}$ to 
$(dm + \bar{m},-d \bar{m},dm' - \bar{m})$.
Let $i \mc M \oplus M \ra M \times I$ be the map that sends 
$(m,m')$ to $(m,0,m')$ and let 
$p \mc M \times I \ra M$ be the map that sends $(m,\bar{m},m')$ to $m+m'$.
One can check easily that $M \times I$ is a chain complex and that
$i$ and $p$ are chain maps whose composite is the codiagonal.
The map $i$ is degreewise split monic with cokernel $\Sigma M$,
so it is a $\Proj$-cofibration, since we have assumed that $M$ is cofibrant.
The map $p$ is a chain homotopy equivalence with chain homotopy inverse
sending $m$ to $(m,0,0)$; 
this implies that it induces a chain homotopy equivalence of generalized 
elements and is thus a $\Proj$-equivalence.
Therefore $M \times I$ is a good cylinder object for $M$.  

It is easy to see that a chain homotopy between two maps $M \ra N$
is the same as a left homotopy using the good cylinder object $M \times I$.
By Lemma~\ref{le:homotopy}, two maps are homotopic if and only if 
they are left homotopic using $M \times I$. 
Thus the model category notion of homotopy is the same
as the notion of chain homotopy when the source is $\Proj$-cofibrant.
\end{proof}

There is a dual proof that proceeds by constructing a specific good path
object $N^{I}$ for $N$ such that a right homotopy using $N^{I}$
is the same as a chain homotopy.

\begin{cor}\label{co:ext}
Let $A$ and $B$ be objects of $\A$ considered as chain complexes
concentrated in degree $0$.
Then $\DA(A,\Sigma^{n} B) \cong \Ext_{\Proj}^{n}(A,B)$.
\end{cor}

See Subsection~\ref{subse:hom-alg} for the definition of the $\Ext$ groups.

\begin{proof}
The group $\DA(A,\Sigma^{n}B)$ may be calculated by choosing a
$\Proj$-cofibrant replacement $A'$ for $A$ and computing the homotopy
classes of maps from $A'$ to $\Sigma^{n} B$.  
(Recall that all objects are $\Proj$-fibrant, so there is no need to take a 
fibrant replacement for $\Sigma^{n} B$.)
A $\Proj$-resolution $P$ of $A$ serves as a 
$\Proj$-cofibrant replacement for $A$, and by Lemma~\ref{le:rh} the homotopy 
relation on $\ChA(P,\Sigma^{n}B)$ is chain homotopy, so it follows that
$\DA(A,\Sigma^{n}B)$ is isomorphic to $\Ext_{\Proj}^{n}(A,B)$.
\end{proof}

More generally, a similar argument shows that the derived functors
of a functor $F$ can be expressed as the cohomology of the derived
functor of $F$ in the model category sense.  To make the story 
complete, we next show that the shift functor $\Sigma$ corresponds
to the notion of suspension that the category $\DA$ obtains as the
homotopy category of a pointed model category.

\begin{defn}\label{de:susp}
Let $\cC$ be a pointed model category.
If $M$ is cofibrant, we define the \dfn{suspension} $\Sigma M$ of $M$ to
be the cofibre of the map $M \coprod M \ra M \times I$ for some
good cylinder object $M \times I$.
(The cofibre of a map $X \ra Y$ is the pushout $* \coprod_X Y$, 
where $*$ is the zero object.)
$\Sigma M$ is cofibrant and well-defined up to homotopy equivalence.
\end{defn}

The loop object $\Omega N$ of a fibrant object $N$ is defined dually.
These operations induce adjoint functors on the homotopy category.
A straightforward argument based on the cylinder object described above
(and a dual path object) proves the following lemma.

\begin{lemma}\label{le:loop}
In the model category $\ChA$,
the functor $\Sigma$ defined in Definition~\ref{de:susp} can be
taken to be the usual suspension, so that
$(\Sigma X)_{n} = X_{n-1}$ and $d_{\Sigma X} = -d_{X}$.
Similarly, $\Omega X$ can be taken to be
the complex $\Sigma^{-1}X$.  That is, $(\Omega X)_{n} = X_{n+1}$
and $d_{\Omega X} = -d_{X}$.  \qed
\end{lemma}

In particular, $\Sigma$ and $\Omega$ are inverse functors.
The second author~\cite{hovey-model} has shown that this implies that cofibre
sequences and fibre sequences agree (up to the usual sign)
and that $\Sigma$ and the cofibre sequences give rise to a 
triangulation of the homotopy category.
(See~\cite[Section I.3]{qu:ha} for the definition of cofibre and
fibre sequences in any pointed model category.)
Using the 
explicit cylinder object from the proof of Lemma~\ref{le:rh},
we can be more explicit.

\begin{cor}\label{co:triangulated}
The category $\DA$ is triangulated with the usual suspension.
A sequence $L \ra M \ra N \ra \Sigma L$ is a triangle if and only if
it is isomorphic in $\DA$ to the usual cofibre sequence on the map $L \ra M$
(see Section~\ref{se:notation}).  \qed
\end{cor}

Now we show that the model structures we construct are proper.
A good reference for proper model categories is~\cite[Chapter~11]{hirschhorn}.

\begin{proposition}\label{prop-proper}
Let $\Proj$ be any projective class on an abelian category $\cA$.
Consider the commutative square in $\Ch{\cA}$ below.
\[
\begin{CD}
A @>f>> X \\
@VqVV @VVpV \\
B @>>g> Y
\end{CD}
\]
\begin{enumerate}
\item [(a)] If the square is a pullback square, $p$ is a
$\Proj$-fibration, and $g$ is a $\Proj$-equivalence, then $f$ is a
$\Proj$-equivalence.  That is, the relative model structure is
\dfn{right proper}.
\item [(b)] If the square is a pushout square,
$q$ is a degreewise split monomorphism, and $f$ is a
$\Proj$-equivalence, then $g$ is a $\Proj$-equivalence.  
In particular, the pushout of a $\Proj$-equivalence along a
$\Proj$-cofibration is a $\Proj$-equivalence.
That is, the relative model structure is \dfn{left proper}.
\end{enumerate}
\end{proposition}

A model category that is both left and right proper is said to be 
\dfn{proper}.
Note that we don't actually need to know that our $\Proj$-cofibrations,
$\Proj$-fibrations and $\Proj$-equivalences give a model structure
to ask whether the structure is proper.

\begin{proof}
Part~(a) is an immediate consequence
of~\cite[Corollary~11.1.3]{hirschhorn}, since every object is
$\Proj$-fibrant.  For part~(b), let $C$ be the cokernel of $q$.
Since pushouts are computed degreewise, it follows that $p$
is a degreewise split monomorphism with cokernel $C$.
Thus we have a map of triangles
\[
\begin{CD}
A @>f>> X \\
@VqVV @VVpV \\
B @>>g> Y \\
@VVV   @VVV \\
C @>>\text{id}> C
\end{CD}
\]
in the homotopy category $\K(\cA)$.  The top and bottom
maps are $\Proj$-equivalences, so the middle map must be as well,
by using the five-lemma and the long exact sequences obtained by
applying the functors $[\Sigma^{k} P, -]$.
\end{proof}

We now consider monoidal structure.  Monoidal model categories are
studied in~\cite[Chapter~4]{hovey-model}. We will assume that $\cA$
is a closed monoidal category.  Thus it is equipped with a functor
$\otimes \mc \cA \times \cA \ra \cA$ such that both $A \otimes -$ and 
$- \otimes A$ have right adjoints for each $A$ in $\cA$.  In particular, 
what we need is that these functors preserve colimits.
There is, of course, an induced closed monoidal structure on $\Ch{\cA}$, 
for which we also use the notation $\otimes$.

\begin{proposition}\label{prop-monoidal}
Let $A$ be a closed monoidal abelian category with a
projective class $\Proj$ such that cofibrant replacements exist and the
unit is $\Proj $-projective.  Then the relative model structure on
$\Ch{\cA}$ is monoidal if and only if the tensor product of two
$\Proj $-cofibrant complexes is always $\Proj $-cofibrant.  
\end{proposition}

\begin{proof}
There are two conditions that must hold for a model category to be
monoidal.  One of them is automatically satisfied when the unit is
cofibrant.  The unit of the monoidal structure on $\Ch{\cA }$ is
$\Sigma^{0} S$, where $S$ is the unit of $\cA$, which has been assumed
to be $\Proj$-projective.  Thus the unit is $\Proj$-cofibrant, by
Lemma~\ref{lem-cofibs}~(a),

Therefore, the relative model structure is
monoidal if and only if whenever $f\mathcolon A\xrightarrow{}B$ and
$g\mathcolon X\xrightarrow{}Y$ are $\Proj $-cofibrations, then the map 
\[
f \boxprod g \mc (A\otimes Y) \amalg _{A\otimes X} (B\otimes X) \lra B\otimes Y
\]
is a $\Proj$-cofibration, and is a $\Proj$-trivial cofibration if either
$f$ or $g$ is a $\Proj$-trivial cofibration.  It is easy to see that
$f\boxprod g$ is a degreewise split monomorphism, with cokernel
$C\otimes Z$.  

By taking $f$ to be the map $0\xrightarrow{}C$ and $g$ to be the map
$0\xrightarrow{}Z$, we see that if the relative model structure is
monoidal, then the tensor product of two $\Proj $-cofibrant complexes is
$\Proj $-cofibrant.  Conversely, if $C\otimes Z$ is $\Proj $-cofibrant
whenever $C$ and $Z$ are $\Proj $-cofibrant, the preceding paragraph
implies $f\boxprod g$ is a $\Proj $-cofibration whenever $f$ and $g$ are
$\Proj $-cofibrations.  If either $f$ or $g$ is a $\Proj $-trivial
cofibration, then one of $C$ or $Z$ is $\Proj $-trivially cofibrant, and
hence contractible.  It follows that $C\otimes Z$ is contractible, and
so $f\boxprod g$ is a $\Proj $-trivial cofibration.  
\end{proof}

In particular, suppose the relative model structure is monoidal, and
$M,N$ are $\Proj $-projective.  Then $\Sigma^{0}M\otimes \Sigma^{0}N\cong
\Sigma^{0}(M\otimes N)$ is $\Proj $-cofibrant, and therefore $M\otimes N$ is
$\Proj $-projective.  With this in mind, we say that a projective class
$\Proj$ on a closed monoidal category $\cA$ is \dfn{monoidal} if the
unit is $\Proj$-projective and the tensor product of $\Proj$-projectives
is $\Proj$-projective.

We do not know of any example of a monoidal projective class $\Proj $
where cofibrant replacements exist, but the relative model structure is
not monoidal.  Certainly this does not happen in either Case A or Case B
of Theorem~\ref{thm-rings}.  In Case A, we have the following
corollary.  

\begin{corollary}\label{cor-monoidal}
Suppose $F \mc \cB \ra \cA$ is a monoidal functor between
closed monoidal abelian categories, with right adjoint $U$
that preserves countable coproducts.
Then the relative model structure on $\Ch{\cA}$ is monoidal.
\end{corollary}

\begin{proof}
Proposition~\ref{prop-mon-case-A} says that $QX\otimes QY$ is $\Proj
$-cofibrant for any $X$ and $Y$, where $Q$ denotes the cofibrant
replacement functor constructed in Section~\ref{se:adjoints}.  If $X$
and $Y$ are already cofibrant, then lifting implies that $X$ is a
retract of $QX$ and $Y$ is a retract of $QY$.  Thus $X\otimes Y$ is a
retract of $QX\otimes QY$, so $X\otimes Y$ is $\Proj $-cofibrant.
Proposition~\ref{prop-monoidal} completes the proof.  
\end{proof}

In Case B, we will show in Corollary~\ref{cor-small-cellular} that every
$\Proj$-cofibrant object is $\Proj$-cellular.  Thus the following
corollary applies.  

\begin{corollary}\label{cor-mon-case-B}
Let $\cA $ be a closed monoidal abelian category with a monoidal
projective class $\Proj $ such that cofibrant replacements exist and
every $\Proj $-cofibrant object is $\Proj $-cellular \ulp
Definition~\ref{de:cellular}\urp .  Then the relative model structure on
$\Ch{\cA }$ is monoidal.  
\end{corollary}

\begin{proof}
Let $A$ and $B$ be $\Proj$-cofibrant.
By Proposition~\ref{prop-monoidal}, it suffices to show that $A \otimes
B$ is $\Proj$-cofibrant.   
By assumption, $A$ is a retract
of a transfinite colimit of a colimit-preserving diagram
\[
0 = A^{0} \lra A^{1} \lra \cdots
\]
such that for each $\alpha$, the map $A^{\alpha} \ra A^{\alpha + 1}$ is
degreewise split monic with cokernel a complex of relative projectives
with zero differential.
Since $- \otimes B$ preserves retracts, degreewise split monomorphisms and
colimits, it is enough to prove that $\Sigma^{i} P \otimes B$ is
$\Proj$-cofibrant for each $\Proj$-projective $P$ and each $i$.
Applying the same filtration argument to $B$, we find that it suffices 
to show that $\Sigma^{i} P \otimes \Sigma^{j} Q = \Sigma^{i+j} (P \otimes Q)$
is $\Proj$-cofibrant for all $\Proj$-projectives $P$ and $Q$ and integers
$i$ and $j$.
But this follows immediately from the fact that $\Proj$ is monoidal
and Lemma~\ref{lem-cofibs}~(a).
\end{proof}

We can also prove a dual statement in Case A.
\comment{Will some form of this go through for any projective class?}

\begin{proposition}\label{prop-injective-monoidal}
Suppose $U \mc \cA \ra \cB$ is a monoidal functor
of closed monoidal abelian categories, with right adjoint $F$.  
Assume that $U$ preserves countable products.
Then the injective relative model structure on $\Ch{\cA}$ is monoidal. 
\end{proposition}

By the ``injective relative model structure'', we mean the model
structure obtained by dualizing Theorem~\ref{thm-rings}.
The (trivial) cofibrations in this model structure are called the
\mdfn{$\cB$-injective (trivial) cofibrations}.

\begin{proof}
Suppose $f\mathcolon A\xrightarrow{}B$ and $g\mathcolon
X\xrightarrow{}Y$ are $\cB$-injective cofibrations with cokernels
$C$ and $Z$, respectively.  This means that $Uf$ and $Ug$ are
degreewise split monomorphisms.  Recall the definition of $f\boxprod
g$ used in the proof of Proposition~\ref{prop-monoidal}.  Since $U$ is
monoidal and preserves 
pushouts, $U(f\boxprod g)\cong Uf\boxprod Ug$.  One can easily check
that $Uf\boxprod Ug$ is a degreewise split monomorphism, so
$f\boxprod g$ is a $\cB$-injective cofibration.  If $f$ is a
$\cB$-injective trivial cofibration, then the cokernel $C$ of $f$
has $UC$ contractible.  Let $Z$ denote the cokernel of $g$.  Since
the cokernel of $Uf\boxprod Ug$ is $UC\otimes UZ$, which is
contractible, $f\boxprod g$ is a $\cB$-injective trivial
cofibration.  Since everything is cofibrant in the injective
relative model structure, the other condition in the definition of a
monoidal model category is automatically satisfied. 
\end{proof}

Another important property of a model category is whether it 
is cofibrantly generated.  This is the topic of 
Section~\ref{se:model-set-small-projs}.

\comment{As far as I am aware, two additional things are needed for a
satisfactory theory of monoids and modules over them.  You need to know
that there is a model category of monoids and a model category of
modules over an arbitrary monoid, and you need to know that a weak
equivalence of monoids gives rise to a Quillen equivalence of the
associated module categories.  In the cofibrantly generated case, the
first problem is solved either by the monoid axiom or by assuming
everything is fibrant and path objects exists for monoids and modules.
When our model categories are cofibrantly generated, everything is
fibrant and path objects certainly exist, so all is well.  But I don't
know what happens in the non-cofibrantly generated case.  For the second
problem, you need to know that, if N is a cofibrant object in the
category of R-modules for a monoid R, then tensoring over R with N takes
weak equivalences of R-modules to weak equivalences.  This does not
appear to depend on cofibrant generation.  The source for all this is
Schwede-Shipley. }

\subsection{Why projective classes?}\label{subse:why}

The astute reader will notice that we haven't used the assumption
that there is always a $\Proj $-epimorphism from a $\Proj$-projective.
Indeed, all we have used is that we have a collection $\Proj$ of
objects such that, with the definitions at the beginning of this
section, cofibrant replacements exist.
The following proposition explains why we start with a projective class.

\begin{prop}
Let $\Proj$ be a collection of objects in an abelian category $\cA$
such that cofibrant replacements exist in $\Ch{\cA}$.  Then there is a unique
projective class $(\Proj', \Epi)$ that gives rise to the
same definitions of weak equivalence, fibration and cofibration.
\end{prop}

\begin{proof}
Let $\Epi$ be the collection of $\Proj$-epimorphisms and
let $\Proj'$ be the collection of all objects $P$ such that
every map in $\Epi$ is $P$-epic.  Then $\Proj'$
contains $\Proj$, and a map is $\Proj'$-epic if and only if
it is $\Proj$-epic.
It follows that the $\Proj'$-exact sequences are the same as the
$\Proj$-exact sequences.
A map $f$ is a $\Proj$-equivalence if and only
if the cofibre of $f$ is $\Proj$-exact.  The same holds
for the $\Proj'$-equivalences, and since the two notions of exactness
also agree, the two notions of equivalence agree.
Finally, since cofibrations are defined in terms of the fibrations
and weak equivalences, the two notions of cofibration agree.

That the pair $(\Proj',\Epi)$ is a projective class follows immediately
from the existence of cofibrant replacements.

Our work earlier in this section shows that the objects of $\Proj'$
are precisely those objects $P$ that are cofibrant when viewed as 
complexes concentrated in degree 0.  Thus the projective class is
the unique projective class giving rise to
the same weak equivalences, fibrations and cofibrations.
\end{proof}

Thus by requiring our collection of objects to be part of a projective
class, we in effect choose a canonical collection of objects determining 
each model structure we produce.  This means that the relevant question is:
which projective classes give rise to model structures?
Theorem~\ref{thm-rings} gives a necessary and sufficient condition,
namely the existence of cofibrant replacements,
and we know of no projective classes that do not satisfy this condition.

An additional advantage of having a projective class is that it
provides us with the language to state results such as:
$\Ext_{\Proj}^{n}(M, N) \cong \Ho\Ch{\cA}(\Sigma^{n} M, N)$
(see Corollary~\ref{co:ext}).

\section{Case A: Projective classes coming from adjoint pairs}
\label{se:adjoints}

In this section we prove Case A of Theorem~\ref{thm-rings}.

Let $U \mc \cA \ra \cB$ be a functor of abelian categories, 
with a left adjoint $F \mc \cB \ra \cA$.  
Then $U$ and $F$ are additive, $U$ is left exact and $F$ is right exact.  
Let $\Proj$ be the projective class on $\cA$ that is the 
pullback of the trivial projective class on $\cB$ 
(see Example~\ref{ex:trivial-projective-class} and 
Subsection~\ref{subse:pullbacks}).
In this subsection we construct a cofibrant replacement functor for the
projective class $\Proj$ under the assumption that $U$ preserves
countable coproducts.

Because this projective class is strong, Lemma~\ref{le:strong}
tells us that the $\Proj$-fibrations and $\Proj$-equivalences
have alternate characterizations.  A map $f$ in $\ChA$ is
a $\Proj$-fibration if and only if $\cA(P,f)$ is degreewise
\emph{split} epic for each $P$ in $\Proj$, 
and is a $\Proj$-equivalence if and only if $\cA(P,f)$ is a 
\emph{chain homotopy equivalence} for each $P$ in $\Proj$.

First we prove a lemma characterizing the fibrations and weak
equivalences in this model structure and giving us a way to generate
cofibrant objects.

\begin{lemma}\label{lem-cofibs-old}
\begin{enumerate}
\item [(a)] A map $p\mathcolon X\xrightarrow{}Y$ is a $\Proj $-fibration
if and only if $Up$ is a degreewise split epimorphism in $\Ch{\cB }$.  
\item [(b)] A map $p\mathcolon X\xrightarrow{}Y$ is a $\Proj
$-equivalence if and only if $Up$ is a chain homotopy equivalence in
$\Ch{\cB }$.  
\item [(c)] If $i$ is a map in $\Ch{\cB}$ that is 
degreewise split monic, then $Fi$ is a $\Proj$-cofibration.
\item [(d)] For any $C$ in $\Ch{\cB}$, $FC$ is $\Proj$-cofibrant.
\end{enumerate}
\end{lemma}

\begin{proof}
For part~(a), note that the $\Proj $-projectives are all retracts of
$FM$ for $M\in \cB $.  Hence $p$ is a $\Proj $-fibration if and only if
$\cA (FM,p)=\cB (M,Up)$ is a surjection for all $M\in \cB $.  This is
true if and only if $Up$ is a degreewise split epimorphism.  For
part~(b), a similar argument shows that $p$ is a $\Proj $-equivalence if
and only if $\cB (M,Up)$ is a quasi-isomorphism for all $M\in \cB $.  We
claim that this forces $Up$ to be a chain homotopy equivalence.  Indeed,
let $C$ denote the cofiber of $Up$.  Then $\cB (M,C)$ is exact for all
$M\in \cB $.  By taking $M=Z_{n}C$, we find that $C$ is exact and that
$C_{n+1}\xrightarrow{}Z_{n}C$ is a split epimorphism.  It follows that
$C$ is contractible, as in the proof of Lemma~\ref{le:strong}.  Thus
$Up$ is a chain homotopy equivalence.  

For part~(c), let $C$ be the cokernel of $i$.  Then $Fi$ is degreewise
split monic with cokernel $FC$.
Suppose that $p\mathcolon X\xrightarrow{}Y$ is a $\Proj$-trivial
fibration with kernel $K$.  Then $K$ is $\Proj $-trivially fibrant, so
part~(b) implies that $UK$ is contractible.
By Lemma~\ref{lem-test-lift}, to show that
$Fi$ is a $\Proj$-cofibration, it suffices to show that
every chain map $FC\xrightarrow{}\Sigma K$ is chain homotopic to $0$.
By adjointness, it suffices to show that every chain map
$C\xrightarrow{}\Sigma UK$ is chain homotopic to $0$.  Since $UK$ is
contractible, this is clear.  Part~(d) follows from part~(c), 
\end{proof}

We now provide a construction, given a complex $X$, of a $\Proj$-cofibrant 
complex $QX$ and a $\Proj$-trivial fibration $QX\xrightarrow{}X$.  This comes
out of the bar construction, which we now recall
from~\cite[Section~IX.6]{maclane-homology}.  Given an object $N$ of
$\cA$, define the complex $BN$ by $(BN)_{m}=(FU)^{m+1}N$ when 
$m \geq -1$ and $0$ otherwise.  We will define maps
$s=s_m \mc U(BN)_{m-1} \ra U(BN)_{m}$ and
$\delta=\delta_m \mc (BN)_m \ra (BN)_{m-1}$.
For $m < 0$ we declare both to be zero.
For $m \geq 0$, $s_m \mc U(BN)_{m-1}=(UF)^mUN \ra U(BN)_m=UF(UF)^mUN$
is defined to be adjoint to the identity map of $F(UF)^mUN$.
For $m \geq 0$, we can then inductively define 
$\delta_{m} \mathcolon (BN)_m=FU(FU)^m N \ra (BN)_{m-1}=(FU)^{m}N$
to be adjoint to the self-map $1-s_{m-1}(U\delta_{m-1})$ of 
$U(FU)^{m}N$.
Properties of adjoint functors then guarantee that 
$(U\delta_{m})s_{m}+s_{m-1}(U\delta_{m-1})=1$.
Using this and the properties of adjoint functors, one deduces a
sequence of implications 
$\delta_{m-1} \delta_{m} = 0 \implies
(U \delta_{m}) s_{m} (U \delta_{m}) = U \delta_{m} \implies
(U \delta_{m})(U \delta_{m+1}) s_{m+1} = 0 \implies 
\delta_{m} \delta_{m+1} = 0$.
Therefore $\delta_{m-1} \delta_{m} = 0$ for each $m$ and so $\delta$
makes $BN$ into a chain complex.

The construction of $BN$, $\delta $, and $s$ is obviously natural in
$N$.  Thus, given a complex $X$ in $\Ch{\cA}$, we get a bicomplex
$(BX)_{m,n}=(BX_{n})_{m}$, where 
\[
\delta \mathcolon (BX)_{m,n}\xrightarrow{}(BX)_{m-1,n} \text{ and }
d\mathcolon (BX)_{m,n}\xrightarrow{}(BX)_{m,n-1}
\]
commute.  Furthermore, $s\mathcolon
U(BX)_{m,n}\xrightarrow{}U(BX)_{m+1,n}$ and $Ud$ also commute.

We can then define a total complex $\overline{Q}X$ by
$(\overline{Q}X)_{k}=\bigoplus _{m+n=k} (BX)_{m,n}$.  The differential
$\partial $ of $\overline{Q}X$ takes the summand $(BX)_{m,n}$ into
$(BX)_{m-1,n}\oplus (BX)_{m,n-1}$ by $(\delta, (-1)^m d)$.  Note that
$\overline{Q}X$ is a filtered complex, with $F^{i}\overline{Q}X$ the
subcomplex consisting of terms with $m\leq i$. In particular,
$F^{-1}\overline{Q}X=\Sigma ^{-1}X$.  We define $QX$ to be
$\overline{Q}X/\Sigma ^{-1}X$, so that $(QX)_{k}=(FU)X_{k}\oplus
(FU)^{2}X_{k-1}\oplus \cdots $.  

The following proposition proves Case A of Theorem~\ref{thm-rings}.

\begin{proposition}\label{prop-bar}
There is a natural $\Proj$-trivial fibration $q_{X}\mathcolon
QX\xrightarrow{}X$, and $QX$ is $\Proj$-cofibrant. 
\end{proposition}

\begin{proof}
We first show that $QX$ is $\Proj$-cofibrant.  There is an increasing
filtration $\{F^{i}QX\}_{i\geq 0}$ on $QX$, where $F^{0}QX=FUX$ and
$F^{i}QX/F^{i-1}QX=(FU)^{i+1}\Sigma ^{i}X$.  Furthermore, each inclusion
$F^{i}QX\xrightarrow{}F^{i+1}QX$ is a degreewise split monomorphism.
By Lemma~\ref{lem-cofibs-old}~(d), each quotient $(FU)^{i+1}\Sigma ^{i}X$ is
$\Proj$-cofibrant, so each map $F^{i-1}QX\xrightarrow{}F^{i}QX$ is a
$\Proj$-cofibration.  Hence $FUX\xrightarrow{}QX$ is a
$\Proj$-cofibration.  Thus $QX$ is $\Proj$-cofibrant.

The map $q_{X}$ is induced by $\delta _{0}\mathcolon
(FU)X_{n}\xrightarrow{}X_{n}$, adjoint to the identity.  The map $q_{X}$
sends the other summands of $(QX)_{n}$ to $0$.  We leave to the reader
the check that this is a chain map.  
Since $U\delta _{0}$ is a split epimorphism, $Uq_{X}$ is a degreewise
split epimorphism, so is a $\Proj$-fibration by
Lemma~\ref{lem-cofibs-old}~(a).  
To show $q_{X}$ is a $\Proj$-equivalence, it suffices to show that the
fiber $\overline{Q}X$ is $\Proj$-contractible, or, equivalently, that
$U\overline{Q}X$ is contractible (Lemma~\ref{lem-cofibs-old}~(b)).  The
contracting homotopy is given by $s$.
Indeed, on the summand $U(BX)_{m,n}$, the
$U(BX)_{m,n}$ component of $s(U\partial) +(U\partial) s$ is $s(U\delta)
+(U\delta) s=1$, and the $U(BX)_{m+1,n-1}$ component is
$s(-1)^{m}(Ud)+(-1)^{m+1}(Ud)s=0$.
(This is where we use that $U$ commutes with coproducts.)
\end{proof}

The following proposition implies that the relative model structure is
monoidal in this case, as explained in Corollary~\ref{cor-monoidal}.  

\begin{proposition}\label{prop-mon-case-A}
Suppose $F\mathcolon \cB \xrightarrow{}\cA $ is a monoidal functor of
closed monoidal abelian categories, with right adjoint $U$ that
preserves countable coproducts.  Then, for any $X$ and $Y$ in $\Ch{\cA
}$, $QX\otimes QY$ is $\Proj $-cofibrant.
\end{proposition}

\begin{proof}
Recall the filtration $F^{i}QX$ on $QX$ used in the proof of
Proposition~\ref{prop-bar}.  Using this filtration, we find that
$QX\otimes QY$ is the colimit of $F^{i}QX\otimes QY$, and each map
$F^{i-1}QX\otimes QY \xrightarrow{}F^{i}QX\otimes QY$ is a degreewise
split monomorphism with cokernel $(FU)^{i+1}\Sigma ^{i}X\otimes QY$.  It
therefore suffcies to show that this cokernel is $\Proj $-cofibrant for
all $i$.  A similar argument using the filtration on $QY$ shows that it
suffices to show that 
\[
(FU)^{i+1}(\Sigma ^{i}X) \otimes (FU)^{j+1}(\Sigma ^{j}Y) \cong
F(U(FU)^{i}\Sigma ^{i}X\otimes U(FU)^{j}\Sigma ^{j}Y) 
\]
is $\Proj $-cofibrant for all $i,j\geq 0$.  But this follow immediately
from Lemma~\ref{lem-cofibs-old}~(d).  
\end{proof}

\subsection{Examples}\label{subse:adjoint-examples}

\begin{example}\label{ex:identity}
Let $\cB$ be a bicomplete abelian category.  Since the identity functor
is adjoint to itself and preserves coproducts, we can apply
Theorem~\ref{thm-rings} to the trivial projective class $\Proj$ to
conclude that $\Ch{\cB}$ is a model category.  We call this the
\mdfn{absolute model structure}.  The $\Proj$-equivalences are the chain
homotopy equivalences (Lemma~\ref{lem-cofibs-old}), the $\Proj$-fibrations
are the degreewise split epimorphisms and the $\Proj$-cofibrations are
the degreewise split monomorphisms.  Every object is both
$\Proj$-cofibrant and $\Proj$-fibrant, and the homotopy category is the
usual homotopy category $\K(\cB)$ in which chain homotopic maps have
been identified.  That this model structure exists was also shown
in~\cite{cole}.  Note that if $\cB $ is closed monoidal, the absolute
model structure is also monoidal, by Corollary~\ref{cor-monoidal}.  In
particular, since every object is $\Proj $-cofibrant, given a
differential graded algebra $R\in \Ch{\cB }$, we get a model structure
on the category of differential graded $R$-modules, where weak
equivalences are chain homotopy equivalences (of the underlying chain
complexes) and fibrations are degreewise split 
epimorphisms.\footnote{This is not obvious; see notes on last page.}

Now let $\cA$ and $\cB$ be bicomplete abelian categories and let 
$U \mc \cA \ra \cB$ be a coproduct preserving functor with left adjoint $F$.
By Theorem~\ref{thm-rings}, $\Ch{\cA}$ has a relative model structure
whose weak equivalences, cofibrations and fibrations are called
$\cB$-equivalences, $\cB$-cofibrations and $\cB$-fibrations.
This structure is a lifting of the absolute model structure on $\Ch{\cB}$
in the sense that a map $f$ in $\Ch{\cA}$ is a weak equivalence
or fibration if and only if $Uf$ is so in $\Ch{\cB}$.
It is often the case that one wants to lift a model structure
along a right adjoint.  
When the model structure is cofibrantly generated, necessary and
sufficient conditions for a lifting are known~\cite[9.1]{dwhika},
~\cite[13.4.2]{hirschhorn}.
Our main theorem says that it is also possible to lift the 
absolute model structure on $\Ch{\cB}$, even though it is not usually 
cofibrantly generated (see Subsection~\ref{subse:not-det-by-a-set}).

It follows from the above that $F$ preserves cofibrations and trivial
cofibrations.  The adjoint functors $F$ and $U$ form a Quillen pair.

The category $\Ch{\cA}$ also has an absolute model structure.
The identity functor sends absolute fibrations and weak equivalences
to $\cB$-fibrations and $\cB$-equivalences.  Thus the identity 
functor is a right Quillen functor from the absolute model structure
to the relative model structure.
\end{example}

\begin{example}\label{ex:change-of-rings-model}
Let $R \ra S$ be a map of rings.  Write \RMod and \SMod for the
categories of left $R$- and $S$-modules.
Consider the forgetful functor $U \mc \SMod \ra \RMod$  and its 
left adjoint $F$ that sends an $R$-module $M$ to $S \tensor_{R} M$.
We saw in Example~\ref{ex:change-of-rings} that this
gives a projective class whose relative projectives are the
$S$-modules $P$ such that the natural map $S \tensor_{R} M \ra M$ 
is split epic as a map of $S$-modules.
The functor $U$ preserves coproducts, so
Theorem~\ref{thm-rings} and Lemma~\ref{le:strong} tell us that 
$\Ch{S} = \Ch{\SMod}$ has
a model structure in which the weak equivalences are the maps
that become chain homotopy equivalences after forgetting the
$S$-module structure.
The fibrations are the maps that in each degree are split epic 
as maps of $R$-modules.
The cofibrations are defined by the \llp the
trivial fibrations.
Equivalently, by Proposition~\ref{prop-cofibrations}, they are the
degreewise split monomorphisms whose cokernels are cofibrant.
And by Lemma~\ref{lem-cofibrant}, a complex $C$ is cofibrant if and only
if each $C_{n}$ is a relative projective and every map from $C$ to
a complex $K$ such that $UK$ is contractible is chain homotopic to $0$.
Lemma~\ref{lem-cofibs} gives us a ready supply of cofibrant objects.
In particular, a relative resolution of an $S$-module $M$ is a
cofibrant replacement, so the group $\Ho\Ch{S}(M, \Sigma^{i} N)$ 
of maps in the homotopy category is isomorphic to the relative $\Ext^{i}$ 
group~\cite{maclane-homology}.

When $R$ and $S$ are commutative, the functor $F$ is monoidal, so this
relative model structure is monoidal by Corollary~\ref{cor-monoidal}.
There is then a derived tensor product $X\otimes ^{L}Y$ in $\Ho
\Ch{S}$, and if $M$ and $N$ are $S$-modules, $H_{i}(M\otimes ^{L}N)$ is
isomorphic to the relative $\Tor _{i}$ group~\cite{maclane-homology}.  
\end{example}

\begin{example}\label{ex:change-of-rings-inj-model}
Again let $R \ra S$ be a map of rings and let $U$ be the forgetful functor.
We saw in Example~\ref{ex:change-of-rings-inj} that $U$ has a right
adjoint $G$ that sends an $R$-module $M$ to $\RMod(S,M)$ and so
we get an injective class on $\SMod$ by pulling back the trivial
injective class along $U$.
$U$ preserves products, so we can apply the duals of Theorem~\ref{thm-rings}
and Lemma~\ref{le:strong} to conclude that $\Ch{S}$ has a model structure
with the same weak equivalences as in the previous example.
The cofibrations are the maps that in each degree are split monic
as maps of $R$-modules, and the fibrations are the maps with the \rlp 
with respect to the trivial cofibrations.
Fibrations and fibrant objects can also be characterized by the duals
of Lemmas~\ref{lem-cofibrant} and~\ref{lem-cofibs}, and the homotopy
category again encodes the relative $\Ext$ groups.  However, this
relative model structure is not monoidal, even when $R$ and $S$ are
commutative, since $U$ is not monoidal.  
\end{example}

\begin{example}\label{ex:bimodules}
Suppose $A$ is an algebra over a commutative ring $k$, and let $\cat{C}$
denote the category of chain complexes of $A$-bimodules.  Then there are
forgetful functors from $A$-bimodules to left $A$-modules, right
$A$-modules, and $k$-modules.  Each of these preserves coproducts, so we
get three different relative model structures on $\cat{C}$.  When we
forget to right or left $A$-modules, the cofibrant replacement functor
$Q$ applied to $A$ gives us the usual un-normalized bar construction,
with $(QA)_{n}=A^{\tensor n+2}$.  The complex $QA$ is then also cofibrant in 
the relative model structure obtained by forgetting to $k$-modules, since it
is a bounded below complex of relative projectives.  The Hochschild
cohomology $HH^{n}(A;N)$ of $A$ with coefficients in a bimodule $N$ (or a
chain complex of bimodules $N$) is then equal to $\Ho\cat{C}(\Sigma^{n}A,N)$, 
where we can use any of the three relative model structures
above.  We can extend this definition by replacing $A$ with an arbitrary
complex of bimodules, but then the answer may depend on which relative
model structure we use.  It is most natural to use the relative model
structure obtained from the forgetful functor to $k$-modules, as then we
can define Hochschild homology as well.

Indeed, the Hochschild homology groups are defined using the tensor
product $M\otimes N$ of bimodules, where we identify $ma\otimes n$ with
$m\otimes an$ but also $am\otimes n$ with $m\otimes na$.  We emphasize
that $M\otimes N$ is only a $k$-module, not a bimodule.   In particular,
this tensor product can't be associative, but it is commutative and also
unital in a weak sense, since $M\otimes (A\otimes A)\cong M$ as a
$k$-module.  Furthermore, the functor $-\otimes N$ has a right adjoint
that takes a $k$-module $L$ to the bimodule $\Hom_{k} (N,L)$, where, if
$g\in \Hom (N,L)$, $(ga)(n)=g(an)$ and $(ag)(n)=g(na)$.  Both of these
functors then extend to functors on chain complexes.  The proof of
Proposition~\ref{prop-monoidal} applies to this case as well, and shows
that if $f$ and $g$ are cofibrations in the relative model structure on
$\cat{C}$ obtained by forgetting to $k$-modules, then $f\boxprod g$ is a
cofibration in the absolute model structure on chain complexes of
$k$-modules.  If either $f$ or $g$ is a trivial cofibration, so is
$f\boxprod g$.  This means the tensor product has a total left derived
functor that is commutative and unital (in the above weak sense).  We
can therefore extend the usual definition of Hochschild homology to
complexes $X$ and $Y$ of bimodules, by defining
$HH_{n}(X,Y)=H_{n}(QX\otimes QY)$.  In fact, it is possible to prove,
using the technique of Proposition~\ref{prop-monoidal}, that, if $X$ is
cofibrant, the functor $X\otimes -$ takes weak equivalences in the
relative model structure to chain homotopy equivalences.  This implies
that $HH_{n}(X,Y)=H_{n}(QX\otimes Y)$.  This is a direct generalization
of the way Hochschild homology is defined
in~\cite[Section~X.4]{maclane-homology}.  
\end{example}

\section{Case B:  Projective classes with enough small projectives}
\label{se:model-small-projs}

In this section we prove that cofibrant replacements exist
in Case B.  We begin by introducing the terminology necessary
for the precise statement of Case B.

We think of an ordinal as the set of all previous ordinals, and
of a cardinal as the first ordinal with that cardinality.

\begin{defn}\label{de:small}
  Given a limit ordinal $\gamma$, the \dfn{cofinality} of $\gamma$,
  $\cofin \gamma $, is the smallest cardinal $\kappa$ such that there exists subset
  $T$ of $\gamma$ of cardinality $\kappa$ with $\sup T = \gamma$.
  The cofinality of a successor ordinal is defined to be $1$.

  A \dfn{colimit-preserving sequence} from an ordinal $\gamma$ to
  a category $\cA$ is a diagram
\[ 
X^{0} \lra X^{1} \lra \cdots \lra X^{\alpha} \lra X^{\alpha + 1} \lra \cdots 
\]
  of objects of $\cA$ indexed by the ordinals less than $\gamma$, such 
  that for each limit ordinal $\lambda$ less than $\gamma$ the natural map 
  $\colim_{\alpha < \lambda} X^{\alpha} \ra X^{\lambda}$ is an isomorphism.

  For $\kappa$ a cardinal, an object $P$ is said to be 
  \mdfn{$\kappa$-small relative to a subcategory $\cM$} if for each 
  ordinal $\gamma$ with $\cofin \gamma > \kappa$ and each colimit-preserving 
  sequence $X \mc \gamma \ra \cA$ that factors through $\cM$, the natural map
  $\colim_{\alpha < \gamma} \cA(P,X^{\alpha}) \ra 
                                \cA(P,\colim_{\alpha < \gamma} X^{\alpha})$
  is an isomorphism, where the last colimit is taken in $\cA$.
\end{defn}

If $\Proj$ is a projective class and $\Proj' \subseteq \Proj$, we
say that the relative projectives in $\Proj'$ are \dfn{enough} if
every $\Proj'$-epimorphism is a $\Proj$-epimorphism.
That is, the collection $\Proj'$ is enough to check $\Proj$-epimorphisms,
$\Proj$-fibrations and $\Proj$-equivalences.
We also sometimes say that $\Proj$ is determined by $\Proj'$.

Fix a cardinal $\kappa$ and let $\Proj'$ be the collection of
all $\Proj$-projectives that are $\kappa$-small relative to
the subcategory of $\cA$ consisting of split monomorphisms with
$\Proj$-projective cokernel.
We say that $\Proj$ has \mdfn{enough $\kappa$-small projectives} if the 
collection $\Proj'$ is enough.

The following proposition proves Case B of Theorem~\ref{thm-rings}.

\begin{prop}\label{pr:enough-small}
Let $\Proj$ be a projective class on a complete and cocomplete
abelian category $\cA$.
Assume that $\Proj$-resolutions can be chosen functorially, and that 
for some cardinal $\kappa$ there are enough $\kappa$-small projectives.
Then functorial cofibrant replacements exist.
\end{prop}

Note that we do not assume that there is a \emph{set} of relative 
projectives that is enough.  The proof below is a variant of the 
small object argument that avoids needing a set of test objects
by using the properties of projective classes.

\begin{proof}
Let $X$ be a chain complex.  The basic construction is the following
``partial cofibrant replacement'':  For each $i$ choose $\Proj$-epimorphisms
$P_i \ra X_i$ and $Q_i \ra Z_i X$ with $P_i$ and $Q_i$ $\Proj$-projective.  
Then form $\bigoplus_i ( D^i P_i \oplus \Sigma^i Q_i )$, which has a natural
map to $X$.  This map is degreewise $\Proj$-epic (because of the $P_i$'s) 
and is epic under $\Ch{\cA}(\Sigma^k P, -)$ for any $\Proj$-projective $P$
(because of the $Q_i$'s).  And its domain is $\Proj$-cofibrant.

We now construct a transfinite sequence $D^1 \ra D^2 \ra \cdots$ with
maps to $X$.  Define $D_1 \ra X$ to be a partial cofibrant replacement.
Let $C^1$ be the pullback
\[
\begin{CD}
C^{1} @>>> D^{1} \\
@VVV @VVV \\
P X @>>> X ,
\end{CD}
\]
where $P X$ is the path complex of $X$.
Let $P^1 \ra C^1$ be a partial cofibrant replacement.  Define
$D^2$ to be the cofibre of the composite $P^1 \ra C^1 \ra D^1$.
The map $D^1 \ra D^2$ is a degreewise split monomorphism whose cokernel 
($\Sigma P^1$) is $\Proj$-cofibrant, so it is a $\Proj$-cofibration.  
Moreover, since the map $P^1 \ra C^1 \ra D^1 \ra X$ factors through $P X$, 
there is a canonical null homotopy, and so there is a canonical map
$D^2 \ra X$.

Let $\gamma$ be an ordinal with cofinality greater than $\kappa$,
and inductively define $D^{\alpha} \ra X$ for $\alpha < \gamma$.
When $\alpha$ is a limit ordinal, $D^{\alpha}$ is defined to be
the colimit of the earlier $D^{\beta}$.
Let $D$ be the colimit of all of the $D^{\alpha}$.
$D$ is $\Proj$-cofibrant and comes with a map to $X$.

Since $D^1 \ra X$ is already degreewise $\Proj$-epic, so is $D \ra X$.
That is, $D \ra X$ is a $\Proj$-fibration.

Let $P$ be a relative projective that is $\kappa$-small relative
to the split monomorphisms with relative projective cokernels.  
We need to show that $D \ra X$ is sent to an isomorphism by 
$[\Sigma^{k} P, -]$.  To simplify notation, we do the case where $k = 0$.
That $[P,D] \ra [P,X]$ is epic is clear since 
$\Ch{\A}(P, D^1) \ra \Ch{\A}(P, X)$ is already so.  
So all that remains is to show that $[P,D] \ra [P,X]$ is monic.  
Let $P \ra D$ be a chain map such that $P \ra D \ra X$ is null.
By Lemma~\ref{le:small-complex} below, the chain map $P \ra D$ factors as a 
chain map through some $D^{\beta}$.  It is null in $X$, so the composite 
$P \ra X$ factors through $P X$.  Thus we get a chain map to the pullback 
$C^{\beta}$.  And this lifts through $P^{\beta} \ra C^{\beta}$.  
So the composite $P \ra D^{\beta} \ra D^{\beta+1}$ must be null.  In
particular, the original map $P \ra D$ must be null.
\end{proof}

We owe the reader a proof of the following lemma. 

\begin{lemma}\label{le:small-complex}
Let $P$ be an object of $\cA$ that is $\kappa$-small relative
to a subcategory $\cM$ of $\cA$.
Then the chain complex $\Sigma^{k} P$ is $\kappa$-small relative
to the subcategory of $\ChA$ consisting of maps whose components
are in $\cM$.
\end{lemma}

\begin{proof}
Suppose $D$ is the colimit of a $\gamma$-indexed diagram whose maps
$D^{\alpha} \ra D^{\alpha+1}$ have components in $\cM$,
where $\gamma$ has cofinality greater than $\kappa$.
For simplicity we treat the case $k = 0$.
Suppose we have a chain map $P \ra D$.
The map $f \mc P \ra D_0$ factors through some $D^{\alpha}_0$, since $P$ 
is small and the maps $D^{\alpha}_{0} \ra D^{\alpha+1}_{0}$ are in $\cM$.
Then $d f \mc P \ra D^{\alpha}_{-1}$ goes to zero in $D$, and so goes to 
zero in some $D^{\beta}$, using the other half of smallness.
So the chain map $P \ra D$ factors as a chain map through some $D^{\beta}$.
This shows that $\colim \ChA(P,D^{\alpha}) \ra \ChA(P,D)$ is surjective.
That it is injective is equivalent to the fact that
$\colim \ChA(P,D^{\alpha}_{0}) \ra \ChA(P,D_{0})$ is injective.
\end{proof}

\begin{corollary}\label{cor-small-cellular}
In the situation of Proposition~\ref{pr:enough-small}, every $\Proj
$-cofibrant complex is $\Proj $-cellular. 
\end{corollary}

\begin{proof}
It suffices to show that the cofibrant replacement $D$ constructed in
Proposition~\ref{pr:enough-small} is $\Proj $-cellular, since if $X$ is
cofibrant, a lifting argument shows that $X$ is a retract of any 
cofibrant replacement of $X$.  
The complex $D$ is constructed as a colimit of a
colimit-preserving functor $D^{\alpha }$, indexed by an ordinal $\gamma$.  
Each map $D^{\alpha }\xrightarrow{}D^{\alpha +1}$ is a
degreewise split monomorphism with cokernel
that is easily seen to be purely $\Proj$-cellular.
By reindexing, one can show that the colimit of such a transfinite diagram
is (purely) $\Proj$-cellular.
\end{proof}

\section{Cofibrant generation}
\label{se:model-set-small-projs}

In this section we assume that $\A$ is a complete and cocomplete
abelian category whose objects are small (Definition~\ref{de:small})
and we characterize the projective classes on $\A$ that give
rise to cofibrantly generated model structures.  We also show 
that for the model structures we put on $\ChA$,
cofibrant generation is equivalent to having a set of weak
generators.  Our goal is the following theorem.

\begin{thm}\label{th:cof-gen-char}
Let $\A$ be a complete and cocomplete abelian category whose
objects are small and let $\Proj$ be any projective class on $\A$.
Then the following are equivalent:
\begin{roenumerate}
\item $\Proj$ is determined by a set %
      \ulp Subsection~\ref{subse:model-cof}\urp .
\item The $\Proj$-equivalences, $\Proj$-fibrations, and $\Proj$-cofibrations
      form a model structure on $\ChA$, and this model structure is
      cofibrantly generated \ulp Subsection~\ref{subse:cof}\urp .
\item The $\Proj$-equivalences, $\Proj$-fibrations, and $\Proj$-cofibrations
      form a model structure on $\ChA$, and the associated homotopy
      category has a set of weak generators 
      \ulp Subsection~\ref{subse:not-det-by-a-set}\urp .
\end{roenumerate}
\end{thm}

\begin{proof}
That (i) $\implies$ (ii) is proved in Subsection~\ref{subse:model-cof}.

That (ii) $\implies$ (iii) is~\cite[Theorem~7.3.1]{hovey-model}.  

That (iii) $\implies$ (i) is proved in Subsection~\ref{subse:not-det-by-a-set}.
\end{proof}

In Subsection~\ref{subse:pure} we give some examples that are
cofibrantly generated, and in Subsection~\ref{subse:not-det-by-a-set} 
we give some examples that are not cofibrantly generated.

\subsection{Background}\label{subse:cof}

In this section we briefly recall the basics of cofibrantly generated
model categories.  This material will be used in the next section to
prove Theorem~\ref{th:cof-gen-char}.  For more details, see the books by Dwyer,
Hirschhorn and Kan~\cite{dwhika}, Hirschhorn~\cite{hirschhorn}, and the
second author~\cite[Section~2.1]{hovey-model}.  We will always assume
our model categories to be complete and cocomplete.  See
Definition~\ref{de:small} for the definition of ``small''.

\begin{defn}
  Let $I$ be a class of maps in a cocomplete category.
  A map is said to be \mdfn{$I$-injective} if it has the right lifting
  property with respect to each map in $I$, and we write $\Iinj$ for the
  category containing these maps.
  A map is said to be an \mdfn{$I$-cofibration} if it has the left lifting
  property with respect to each map in $\Iinj$, and we write $\Icof$ for the
  category containing these maps.
  A map is said to be \mdfn{$I$-cellular} if it is a transfinite composite
  of pushouts of coproducts of maps in $I$, and we write $\Icell$ for the
  category containing these maps.
\end{defn}

Note that $\Icell$ is a subcategory of $\Icof$.

\begin{defn}
  A \dfn{cofibrantly generated model category} is a model category
  $\cM$ for which there exist sets $I$ and $J$ of morphisms with
  domains that are small relative to $\Icof$ and $\Jcof$, respectively,
  such that $\Icof$ is the category of cofibrations and
  $\Jcof$ is the category of trivial cofibrations.
  It follows that $\Iinj$ is the category of trivial fibrations
  and that $\Jinj$ is the category of fibrations.
\end{defn}

For example, take $\cM$ to be the category of spaces and take
$I = \{S^{n}\ra B^{n+1}\}$ and $J = \{B^{n} \times 0 \ra B^{n} \times [0,1]\}$.

\begin{prop} \textup{(Recognition Lemma.)}\label{pr:rl}
  Let $\cM$ be a category that is complete and cocomplete, let
  $W$ be a class of maps that is closed under retracts and
  satisfies the two-out-of-three axiom, and let $I$ and $J$ be
  sets of maps with domains that are small relative to 
  $\Icell$ and $\Jcell$, respectively, such that
  \begin{roenumerate}
  \item $\Jcell \subseteq \Icof \cap W$ and $\Iinj \subseteq \Jinj \cap W$, and
  \item $\Jcof \supseteq \Icell \cap W$ or $\Iinj \supseteq \Jinj \cap W$.
  \end{roenumerate}
  Then $\cM$ is cofibrantly generated by $I$ and $J$, and $W$ is the
  subcategory of weak equivalences. 
  Moreover, a map is in $\Icof$ if and only if it is a retract
  of a map in $\Icell$, and similarly for $\Jcof$ and $\Jcell$.
  \qed
\end{prop}

The proof, which uses the small object argument and is due to Kan, can
be found in~\cite{dwhika}, \cite{hirschhorn},
and~\cite[Theorem~2.1.19]{hovey-model}.  

\subsection{Projective classes with sets of enough small 
projectives}\label{subse:model-cof}

In this subsection we prove that under a hypothesis slightly
stronger than Case B of Theorem~\ref{thm-rings} we can conclude
that the model structure is cofibrantly generated.
We will prove this from scratch, since the additional work is small,
thanks to the theory of cofibrantly generated model categories.

Let $\A$ be a complete and cocomplete abelian category
equipped with a projective class $\Proj$.
We use the terminology from Section~\ref{se:notation} and
the refer the reader to Definition~\ref{de:P-structure} for the
definitions of $\Proj$-equivalence, $\Proj$-fibration and
$\Proj$-cofibration.

The following lemma will be used to prove our main result.

\begin{lemma}\label{le:injs}
Consider an object $P \in \A$ and a map $f \mc X \ra Y$ of chain complexes.
\begin{roenumerate}
\item The map $f$ has the \rlp each of the
maps $0 \ra D^{n}P$, $n \in \Z$, if and only if 
the induced map of $P$-elements is a surjection.
\item The map $f$ has the \rlp each of the
maps $\Sigma^{n-1} P \ra D^{n}P$, $n \in \Z$, if and only if 
the induced map of $P$-elements is a surjection and a quasi-isomorphism.
\end{roenumerate}
\end{lemma}

We use the following terminology in the proof.
The cycles, boundaries and homology classes in
$\A(P,X)$ are called \mdfn{$P$-cycles}, \mdfn{$P$-boundaries} and
\mdfn{$P$-homology classes} in $X$.  Note that a $P$-cycle is the same
as a $P$-element of $Z_{n} X$, but that not every $P$-element of
$B_{n} X$ is necessarily a $P$-boundary.
A $P$-element is just a chain map $D^{n} P \ra X$,
a $P$-cycle is just a chain map $\Sigma^{n} P \ra X$, with $P$
regarded as a complex concentrated in degree $0$, and a
$P$-homology class is just a chain homotopy class of maps $\Sigma^{n} P \ra X$.

\begin{proof}
We begin with (i): The map $f$ has the \rlp the
map $0 \ra D^{n}P$ if and only if each map $D^{n}P \ra Y$ factors
through $f$, \ie if and only if each $P$-element of $Y_{n}$
is in the image of $f_{n}$.

Now (ii): The map $f$ has the \rlp 
the map $\Sigma^{n-1}P \ra D^{n}P$ if and only if 
for each $P$-element $y$ of $Y_{n}$ whose boundary is the image of
a $P$-cycle $x$ of $X_{n-1}$, there is a $P$-element $x'$ of 
$X_{n}$ that hits $y$ under $f$ and $x$ under the differential.
(In other words, if and only if the natural map 
$X_{n} \ra Z_{n-1}X \times_{Z_{n-1}Y} Y_{n}$
induces a surjection of $P$-elements.)

So suppose that $f$ has the \rlp each map
$\Sigma^{n-1}P \ra D^{n}P$.
As a preliminary result, 
we prove that $f$ induces a surjection of $P$-cycles.
Suppose we are given a $P$-cycle $y$ of $Y$.
Its boundary is zero and is thus the image of the $P$-cycle $0$ of $X$.
Therefore $y$ is the image of a $P$-cycle $x'$.

It follows immediately that $f$ induces a surjection in $P$-homology.

Now we prove that $f$ induces a surjection of $P$-elements.
Suppose we are given a $P$-element $y$ of $Y$.
By the above argument, its boundary, which is a $P$-cycle,
is the image of a $P$-cycle $x$ of $X$.
Thus, by the characterization of maps $f$ having the RLP,
we see that there is a $P$-element $x'$ that hits $y$.

A similar argument shows that $f$ induces an injection in $P$-homology.

We have proved that if $f$ has the \rlp the maps
$\Sigma^{n-1}P \ra D^{n}P$, then $f$ induces an isomorphism in $P$-homology
and a surjection of $P$-elements.
The proof of the converse goes along the same lines.
\end{proof}

\begin{cor}\label{co:char}
A map $f \mc X \ra Y$ is a $\Proj$-fibration 
\ulp resp.\ $\Proj$-trivial fibration\urp\ 
if and only if it has the \rlp the map $0 \ra D^{n}P$ 
\ulp resp.\ $\Sigma^{n-1} P \ra D^{n} P$\urp\ 
for each $\Proj$-projective $P$ and each $n \in \Z$.    \qed
\end{cor}

We want to claim that the above definitions lead to a cofibrantly generated
model category structure on the category $\ChA$.  In order to
prove this, we need to assume that there is a \emph{set} $\cS$
of $\Proj$-projectives such that a map $f \mc A \ra B$
is $\Proj$-epic if and only if 
$f$ induces a surjection of $P$-elements for each $P$ in $\cS$.
This implies that for any $B$ there is a $\Proj$-epimorphism $P \ra B$
with $P$ a coproduct of objects from $\cS$, and that every 
$\Proj$-projective object is a retract of such a coproduct
(cf.\ Lemma~\ref{le:det-by-a-set}).
We also need to assume that each $P$ in $\cS$ is small relative to the 
subcategory $K$ of split monomorphisms with $\Proj$-projective cokernels.
(See Section~\ref{se:model-small-projs} for the definition of ``small.'')
When these conditions hold we say that $\Proj$ is 
\dfn{determined by a set of small objects} or that there
is a \dfn{set of enough small projectives}.

\begin{thm}\label{th:dc}
Assume that $\Proj$ is determined by a set $\cS$ of small objects.
Then the category $\ChA$ is a cofibrantly generated model category
with the following generating sets:
\[ I := \{ \Sigma^{n-1} P \ra D^{n} P \st P \in \cS,\, n \in \Z \} , \]
\[ J := \{              0 \ra D^{n} P \st P \in \cS,\, n \in \Z \} . \]
The weak equivalences, fibrations and cofibrations are the
$\Proj$-equivalences, $\Proj$-fibrations and $\Proj$-cofibrations
as described in Definition~\ref{de:P-structure}.
A map is a $\Proj$-cofibration if and only if it is degreewise split monic
with $\Proj$-cellular cokernel.
\end{thm}

In particular, every object is $\Proj$-fibrant and an object is 
$\Proj$-cofibrant if and only if it is $\Proj$-cellular 
(Definition~\ref{de:cellular}).

\begin{proof}
We check the hypotheses of the Recognition Lemma from the previous subsection.
Since $\cA$ is complete and cocomplete, so is $\ChA$;  limits and
colimits are taken degreewise.
The class $W$ of $\Proj$-equivalences is easily seen to be closed under 
retracts and to satisfy the two-out-of-three condition.
The zero chain complex is certainly small relative to $\Jcell$.
It is easy to see that a map is in $\Icell$ if and only if it is a 
degreewise split monomorphism whose cokernel is purely cellular
(Definition~\ref{de:cellular}).
In particular, every map in $\Icell$ has components in $K \subseteq \cA$.
We assumed that each $P$ in $\cS$ is small relative to $K$, 
and so it follows from Lemma~\ref{le:small-complex} that each $\Sigma^{k} P$ 
is small relative to $\Icell \subseteq \ChA$.

Since the projectives in $\cS$ are enough to test whether
a map is a $\Proj$-fibration or a $\Proj$-equivalence, Corollary~\ref{co:char}
tells us that $\Iinj$ is the collection of $\Proj$-trivial fibrations
and that $\Jinj$ is the collection of $\Proj$-fibrations.
Thus we have an equality $\Iinj = \Jinj \cap W$, giving us two
of the inclusions required by the Recognition Lemma.

We now prove that $\Jcell \subseteq \Icof \cap W$.
Since $\Iinj \subseteq \Jinj$, it is clear that $\Jcof \subseteq \Icof$,
so in particular $\Jcell \subseteq \Icof$.
We must prove that $\Jcell \subseteq W$, \ie that each map that
is a transfinite composite of pushouts of coproducts of
maps in $J$ is a $\Proj$-equivalence.
A map in $J$ is of the form $0 \ra D^{n} P$ for some $P \in \cS$.
Thus a pushout of a coproduct of maps in $J$ is of the form 
$X \ra X \oplus C$, with $C$ a contractible complex, and
a transfinite composite of such maps is of the same form as well.
Thus such a map is in fact a chain homotopy equivalence, so it is
certainly a $\Proj$-equivalence.

We can now apply the Recognition Lemma and conclude that $\ChA$
is a model category with weak equivalences the $\Proj$-equivalences.
The fibrations are the maps in $\Iinj$, which, as we noted above,
are the $\Proj$-fibrations.
The cofibrations are the maps in $\Icof$, \ie the maps with the
\llp the $\Proj$-trivial fibrations, and this is precisely
how the $\Proj$-cofibrations were defined.

The recognition lemma tells us that the $\Proj$-cofibrations consist 
precisely of the retracts of maps in $\Icell$.
As discussed above, $\Icell$ is the class of degreewise split monomorphisms
with purely cellular cokernels.
So the $\Proj$-cofibrations are the degreewise split monomorphisms 
whose cokernels are cellular.
\end{proof}

\subsection{The pure and categorical derived categories}\label{subse:pure}

In this section we let $R$ be an associative ring with unit and 
we take for $\A$ the category of left $R$-modules.
We are concerned with two projective classes on the category $\A$.
The first is the \dfn{categorical projective class} $\cC$ 
whose projectives are summands of free modules, 
whose exact sequences are the usual exact sequences,
and whose epimorphisms are the surjections.
The second is the \dfn{pure projective class} $\cP$ 
whose projectives are summands of sums of finitely presented modules.
A short exact sequence is $\Proj$-exact iff it remains exact
after tensoring with any right module.
A map is a $\Proj$-epimorphism iff it appears as the epimorphism in a $\Proj$-exact
short exact sequence.
We say \dfn{pure projective} instead of $\Proj$-projective, and similarly
for \dfn{pure exact} and \dfn{pure epimorphism}.
We assume that the reader has some familiarity with these
projective classes.  A brief summary with further references
may be found in~\cite[Section~9]{ch:itcpgs}.  As usual, we write
$\Ext^{*}(-,B)$ (resp.\ $\PExt^{*}(-,B)$) for the derived functors of 
$\A(-,B)$ with respect to the categorical (resp.\ pure) projective 
class. 

Both of these projective classes are determined by sets of small
objects:  $\cC$ is determined by $\{R\}$ and $\cP$ is determined
by any set of finitely presented modules containing a representative
from each isomorphism class.
Thus we get two cofibrantly generated model category structures on $\ChR$ 
and two derived categories, the categorical derived category $\DC$ and
the pure derived category $\DP$, both containing a set of weak generators.
We refer to the pure weak equivalences in $\ChR$ as \dfn{pure quasi-isomorphisms},
and as usual call the categorical weak equivalences simply quasi-isomorphisms.
Similarly, we talk of pure fibrations and fibrations, pure cofibrations
and cofibrations, etc.

Pure homological algebra is of interest in stable homotopy
theory because of the following result.

\begin{thm}\cite{chst:pmht}
Phantom maps from a spectrum $X$ to an Eilenberg-Mac\,Lane spectrum $HG$
are given by $\PExt_{\Z}^{1}(H_{-1} X, G)$, that is, by maps of degree
one from $H_{-1} X$ to $G$ in the pure derived category of abelian groups.
\end{thm}

In addition to the connection between phantom maps and pure homological
algebra, the authors are interested in the pure derived category as a
tool for connecting the global pure dimension of a ring $R$ to
the behaviour of phantom maps in $\DC$ and $\DP$ under composition.

The two derived categories are connected in various ways.
For example, since every categorical projective is pure projective, it follows
that every pure quasi-iso\-mor\-phism is a quasi-isomorphism.  This
implies that there is a unique functor $R: \DP \ra \DC$
commuting with the functors from $\ChR$.  This functor has a left
adjoint $L$.  We can see this in the following way:  the identity
functor on $\ChR$ is adjoint to itself.  It is easy to see that
every pure fibration is a fibration and that every pure trivial
fibration is a trivial fibration.  Therefore, 
by~\cite[Theorem~9.7]{dwsp:htmc}, there is an induced pair
of adjoint functors between the pure derived category and
the categorical derived category, and the right adjoint is the functor $R$
mentioned above. 
This right adjoint is the identity on objects, since everything is
fibrant, and induces the natural map $\PExt^{*}(A,B) \ra \Ext^{*}(A,B)$
for $R$-modules $A$ and $B$.
The left adjoint sends a complex $X$ in $\DC$ to
a (categorical) cofibrant replacement $\tilde{X}$ for $X$.
Similar adjoint functors exist whenever one has two projective
classes on a category, one containing the other.

\subsection{Failure to be cofibrantly generated}\label{subse:not-det-by-a-set}

In this subsection, we show that many of the model structures we have
constructed are not cofibrantly generated and that their homotopy
categories do not have a set of weak generators (see below).
In particular, we show this for the absolute model structures on
$\Ch{\Z}$ and $\Ch{\Z_{(p)}}$.
Recall that the weak equivalences in these model structures are the
chain homotopy equivalences, the cofibrations are the degreewise split
monomorphisms, and the fibrations are the degreewise split
epimorphisms.

A set $\cG$ of objects in an additive category $\cH$ is a \dfn{set of
weak generators} if each non-zero object $X$ in $\cH$ receives a
non-zero map from some $G$ in $\cG$.

A number of people have recently found proofs that certain model
categories are not cofibrantly generated.  The first such proof we
have heard of is due to Dan Isaksen (personal communication), who
proved that his model structure on pro-simplicial sets~\cite{isaksen} is
not cofibrantly generated.  In addition, Ad\'{a}mek, Herrlich,
Rosick\'{y}, and Tholen~\cite{adamek-et-al} have constructed another
non-cofibrantly generated model structure.  Furthermore, Neeman has
proved~\cite[Lemma~E.3.2]{ne:tc} that the triangulated category
$\cat{K}(\Z )$ does not have a set of weak generators.  This implies
that the absolute model structure on $\Ch{\Z }$ is not cofibrantly
generated, by~\cite[Theorem~7.3.1]{hovey-model}.

We also prove that our model structures are not cofibrantly generated
by showing that their homotopy categories do not have sets of weak
generators.  And we prove the latter using the following result, which
completes the proof of Theorem~\ref{th:cof-gen-char}.

\begin{thm}\label{th:weak-gens}
Let $\A$ be a complete and cocomplete abelian category whose objects
are small and let $\Proj$ be any projective class on $\A$.  If the
$\Proj$-equivalences, $\Proj$-fibrations, and $\Proj$-cofibrations
form a model structure on $\ChA$, and the associated homotopy category
has a set of weak generators, then $\Proj$ is determined by a set of
small objects.
\end{thm}

For the term ``small,'' see Definition~\ref{de:small}.
In our examples, $\A$ will be the category $\RMod$ for some ring $R$,
and in this category every object is small.

\begin{proof}
Let $\cG$ be a set of weak generators for $\Ho\ChA$ and, without loss
of generality, assume that each $G$ in $\cG$ is $\Proj$-cofibrant.
Let $\cS = \{ G_{n} \st n \in \Z \text{ and } G \in \cG \}$.
Then each $S$ in $\cS$ is $\Proj$-projective.
Suppose $p \mc X \ra Y$ in $\A$ is $\cS$-epic.  
We must show that $p$ is $\Proj$-epic.
Consider the chain complex $\ker p \ra X \ra Y$, with zeroes
elsewhere.
Because the map $p$ is $\cS$-epic, it is easy to check that
every map from a generator $G$ to the complex $\ker p \ra X \ra Y$
is null homotopic.
But since $\cG$ is a set of weak generators, this implies that
this complex is $\Proj$-equivalent to the zero complex.
That is, the complex $0 \ra \A(P,\ker p) \ra \A(P,X) \ra \A(P,Y) \ra 0$
is exact for each $P$ in $\Proj$.  In particular, the map $X \ra Y$ is 
$\Proj$-epic.
\end{proof}

We now give examples of projective classes that are not
determined by a set.
We begin with a lemma about abelian groups. 

\begin{lemma}\label{lem-non-cof}
For any cardinal $\kappa $, there is an abelian group $A$ such that $A$
is not a retract of any direct sum of abelian groups of cardinality less
than $\kappa $.  
\end{lemma}

\begin{proof}
We use the Ulm invariants, as described in~\cite{kaplansky}.  Recall
that for each prime $p$, each abelian group $A$ and each ordinal
$\lambda $, the Ulm invariant $U_{p}(A,\lambda )$ is a cardinal number.
 From the definition it is clear that $U_{p}(A,\lambda )=0$ for $\lambda
$ larger than the cardinality of $A$.  Furthermore, $U_{p}(-,\lambda )$
takes direct sums of abelian groups to sums of
cardinals~\cite[Problem~31]{kaplansky}.  Hence we only need to find an
abelian group $A$ with $U_{p}(A,\kappa )\neq 0$.  Such a ($p$-torsion)
group exists for every $p$ by~\cite[Problem 43]{kaplansky}.
\end{proof}

Recall now the trivial projective class $\Proj$ on a category $\A$,
which has all objects $\Proj$-projective and only the split epimorphisms
$\Proj$-epic.

\begin{lemma}\label{lem-non-cof-maps}
The trivial projective class on the category of abelian groups
is not determined by a set.
\end{lemma}

\begin{proof}
For any set $\cS$ of abelian groups, we must exhibit a map $p \mc X \ra Y$
that is $\cS$-epic but not split epic.
So fix a set $\cS$ and let $\kappa$ be a cardinal larger than the 
cardinality of each group in $\cS$.
Let $Y$ be an abelian group that is not a retract of any direct sum of
abelian groups of cardinality less than $\kappa $, using
Lemma~\ref{lem-non-cof}.  Let $X$ denote the direct sum of
the $C$ in $\cS$, with one copy of $C$ for each homomorphism
$C\xrightarrow{}Y$.  Then there is an obvious map $p\mathcolon
X\xrightarrow{}Y$, and this map is not split epic since $Y$ is not a
retract of $X$.  However, if $B$ is an abelian group in $\cS$,
any homomorphism $f\mathcolon B\xrightarrow{}Y$
obviously factors through $p$.  
\end{proof}

 From this lemma and Theorem~\ref{th:cof-gen-char} we deduce:

\begin{cor}\label{co:ChZ-non-cof}
The absolute model structure on $\Ch{\Z}$ is not cofibrantly
generated, and the homotopy category $\cat{K}(\Z)$ does
not have a set of weak generators.  \qed
\end{cor}

The results above generalize.  Call a ring $R$ \dfn{difficult} if 
the trivial projective class on the category of $R$-modules is
not determined by a set.
We showed above that $\Z$ is difficult, and it seems to be
the case that many rings are difficult.
For example, every proper subring $R$ of $\Q$ is difficult.  
Indeed, the proof of Lemma~\ref{lem-non-cof} goes through with $p$ 
chosen to be a prime that is not invertible in $R$.
Then the proof of Lemma~\ref{lem-non-cof-maps} applies to show that
$R$ is difficult.
Similarly, one can show that polynomial rings are difficult,
by using an indeterminate in place of the prime $p$.
\comment{In fact, I think the whole Ulm argument, and the
Problems of Kaplansky, go through for any ring with a regular
element (non-unit, non-zero divisor), but I have to check the
details. }%
Moreover, we have the following lemma.

\begin{lemma}\label{le:difficult-rings}
If $R$ is a difficult ring and $R \ra S$ is a map of rings that
is split monic as a map of $R$-modules, then $S$ is difficult.
In addition, the \emph{relative} projective class on $\SMod$
is not determined by a set.
\end{lemma}

\begin{proof}
Fix a cardinal $\kappa$.  Since $R$ is difficult, there is a map $p$ of
$R$-modules that is not split epic such that $\Hom(B,p)$ is epic
whenever $B$ has cardinality less than $\kappa$.  Let $q$ be the map
$\Hom_{R}(S,p)$ of $S$-modules.  Since $R$ splits off of $S$, $q$ is the
sum of $p$ and another map as a map of $R$-modules.  Thus $q$ is not
split epic as an $R$-module map and hence as an $S$-module map.  Since
$\Hom _{R}(S,-)$ is right adjoint to the forgetful functor, we have
$\Hom _{S}(M,q)=\Hom _{R}(M,p)$ for any $S$-module $M$.  In particular,
$\Hom _{S}(M,q)$ is surjective for all $S$-modules $M$ of cardinality
less than $\kappa $.  Since any set of $S$-modules will have sizes
bounded above by some cardinal, we have shown that the ring $S$ is
difficult.

Next we show that the relative projective class $\Proj$ on $\SMod$ is 
also not determined by a set.  Recall that the $\Proj$-projectives
are the retracts of extended modules $M \tensor_{R} S$
and the $\Proj$-epimorphisms are the maps of $S$-modules that
are split epic as $R$-module maps.
So for each $\kappa$ we must exhibit a map $q$ of $S$-modules
that is not split epic as an $R$-module map such that, for every extended
module $M \tensor_{R} S$ of cardinality less than $\kappa$,
the map $\Hom_{S}(M \tensor_{R} S, q)$ is epic.
Note that neither requirement uses the $S$-module structure on $q$.
We choose $q$ as in the previous paragraph, where we already noted
that it is not split epic as an $R$-module map.  Since $\Hom _{S}(N,q)$ is
surjective for all $S$-modules $N$ of cardinality less than $\kappa $,
of course $\Hom _{S}(M\tensor _{R}S,q)$ is surjective for all extended
modules $M\tensor _{R}S$ of cardinality less than $\kappa $.  
\end{proof}

On the other hand, if $R$ is a semisimple ring, or equivalently, 
if every $R$-module is projective, then the trivial projective
class is the same as the categorical projective class, and
so is determined by a set (the set $\{R\}$).
It follows that the absolute model structure on $\Ch{R}$ 
is cofibrantly generated and that $\KR = \DR$ has a set of 
weak generators.
The hypotheses hold for fields and the rings $\Z/n$, where
$n$ is a product of distinct primes, for example.
For other $n$, the ring $\Z/n$ is not semisimple, but it is
easy to see that the trivial projective class is again
determined by a set.

\comment{Can we characterize the difficult rings?  Are they just the
non-semisimple rings?  No.  Any other guesses?}

\section{Simplicial objects and the bounded below derived category}\label{se:simp}

In this section, we discuss the notion of a projective class on
a pointed category and prove that under certain conditions the
category of simplicial objects has a model structure reflecting
a given projective class.  We use this to deduce that for any
projective class on any complete and cocomplete abelian category,
the category of bounded below chain complexes has a model structure
reflecting the projective class.  We do not need to assume that
our projective class comes from an adjoint pair or that there are 
enough small projectives.  We also deduce that the category of
$G$-equivariant simplicial sets has a model structure for any
group $G$ and any family of subgroups.

\subsection{Projective classes in pointed categories}\label{subse:pc-ptd}

One can define the notion of a projective class on any pointed
category by replacing the use of epimorphisms with exact sequences.
However, in our applications we will assume that our pointed categories
are bicomplete;  in particular, they will have kernels, and when this
is in the case, the definition given in Section~\ref{se:pc} is
equivalent to the more general definition~\cite{eimo:frha}.
So in what follows, a projective class on a pointed category
is a collection $\Proj$ of objects and a collection $\Epi$ of maps
satisfying Definition~\ref{de:pc}.

We supplement the examples from Section~\ref{se:pc} with the
following non-additive examples.

\begin{example}\label{ex:pointed-sets}
Let $\A$ be the category of pointed sets, let $\Proj$ be the collection
of all pointed sets and let $\Epi$ be the collection of all epimorphisms.
Then $\Epi$ consists of the $\Proj$-epimorphisms and
these classes form a projective class.
\end{example}

\begin{example}\label{ex:trivial-projective-class-ptd}
More generally, if $\A$ is any pointed category with kernels, 
$\Proj$ is the collection of all objects, 
and $\Epi$ is the collection of all split epimorphisms,
then $\Proj$ is a projective class.  
It is called the \dfn{trivial projective class}.
\end{example}

Example~\ref{ex:pointed-sets} %
is determined by a set in the sense of Lemma~\ref{le:det-by-a-set}.

As in the abelian case, 
a projective class is precisely the information needed to form
projective resolutions and define derived functors.
However, in the non-additive case not all of the usual
results hold.

For further details we refer the reader to the
classic reference~\cite{eimo:frha}.

\subsection{The model structure}

Let $\cA$ be a complete and cocomplete pointed category
with a projective class $\Proj$.  We do not assume that $\cA$ is abelian.
Write $\sA$ for the category of simplicial objects in $\cA$.
Given a simplicial object $X$ and an object $P$, write $\cA(P,X)$
for the simplicial set that has $\cA(P,X_{n})$ in degree $n$.
Define a map $f$ to be a \mdfn{$\Proj$-equivalence} 
(resp.\ \mdfn{$\Proj$-fibration})
if $\cA(P,f)$ is a weak equivalence (resp.\ fibration) of simplicial sets 
for each $P$ in $\Proj$.
Define $f$ to be a \mdfn{$\Proj$-cofibration} if it has the \llp 
the $\Proj$-trivial fibrations.

Consider the following conditions:
\begin{itemize}
\item [(*)] For each $X$ in $\sA$ and each $P$ in $\Proj$, $\cA(P,X)$
is a fibrant simplicial set.
\item [(**)] $\Proj$ is determined by a set $\cS$ of small objects.
Here we require that each $P$ in $\cS$ be small\footnote{We in fact 
need that the objects in $\cS$ are $\omega$-small; see last page.}
\enlargethispage{8pt}
with respect to all
split monomorphisms in $\cA$, not just the split monomorphisms with
$\Proj$-projective cokernels.
\end{itemize}

\begin{thm}\label{th:simp}
Let $\cA$ be a complete and cocomplete pointed category with a 
projective class $\Proj$.  If $\cA$ and $\Proj$ satisfy \ulp *\urp\ or 
\ulp **\urp, then the $\Proj$-equivalences, $\Proj$-fibrations and 
$\Proj$-cofibrations form a simplicial model structure.
When \ulp **\urp\ holds, this model structure is cofibrantly generated 
by the sets
\[ I := \{ P \tensor \dot{\Delta}[n] \ra P \tensor \Delta[n] 
       \st P \in \cS,\, n \geq 0 \} , \]
\[ J := \{ P \tensor V[n,k] \ra P \tensor \Delta[n] 
       \st P \in \cS,\, 0 < n \geq k \geq 0 \} . \]
\end{thm}

In the above, $\Delta[n]$ denotes the standard $n$-simplex simplicial set;  
$\dot{\Delta}[n]$ denotes its boundary;
and $V[n,k]$ denotes the subcomplex with the $n$-cell and its
$k$th face removed.
For an object $P$ of $\cA$ and a set $K$, $P \tensor K$ denotes
the coproduct of copies of $P$ indexed by $K$.
When $K$ is a simplicial set, $P \tensor K$ denotes the simplicial
object in $\cA$ that has $P \tensor K_{n}$ in degree $n$.

If $\cA$ is abelian, then (*) holds.  Moreover, in this case
we have the Dold-Kan equivalence given by
the normalization functor $\sA \ra \textup{Ch}^{+}(\cA)$,
where $\textup{Ch}^{+}(\cA)$ denotes the category of non-negatively graded
chain complexes of objects of $\cA$.  
Thus we can deduce:

\begin{cor}
Let $\cA$ be a complete and cocomplete abelian category with
a projective class $\Proj$.  Then $\textup{Ch}^{+}(\cA)$ is a model category.
A map $f$ is a weak equivalence iff $\cA(P,f)$ is
a quasi-isomorphism for each $P$ in $\Proj$.
A map $f$ is a fibration iff $\cA(P,f)$ is surjective in positive
degrees \ulp but not necessarily in degree 0\urp\ for each $P$ in $\Proj$.
A map is a cofibration iff it is degreewise split monic with degreewise
$\Proj$-projective cokernel. 
Every complex is fibrant, and a complex is cofibrant iff
it is a complex of $\Proj$-projectives. \qed
\end{cor}

Note that no conditions on the projective class are required, and that
the description of the cofibrations is simpler than in the unbounded case.

Special cases of the theorem and corollary were proved by 
Blanc~\cite{bl:aiht}, with the hypothesis that the $\Proj$-projectives
are cogroup objects.  
(Under this hypothesis, condition (*) automatically holds.)

As another example of the theorem, let $G$ be a group and consider
pointed $G$-simplicial sets, or equivalently, simplicial objects
in the category $\cA$ of pointed (say, left) $G$-sets.
Let $\cF$ be a family of subgroups of $G$, and consider the
set $\cS = \{ (G/H)_{+} \st H \in \cF \}$ of homogeneous spaces
with disjoint basepoints.
These are small, and determine a projective class.
Thus, using case (**) of the above theorem, we can deduce:

\begin{cor}\label{co:G-sSet}
Let $G$ be a group\footnote{$G$ should be finite; see last page.} 
\enlargethispage{4pt}
and let $\cF$ be a family of subgroups of $G$.
The category of pointed $G$-equivariant simplicial sets has a model
category structure in which the weak equivalences are precisely
the maps that induce a weak equivalence on $H$-fixed points
for each $H$ in $\cF$.  \qed
\end{cor}

We omit the proof of Theorem~\ref{th:simp}, and simply note that
it follows the argument in \cite[Section II.4]{qu:ha} fairly
closely.  It is a bit simpler, in that Quillen spends
part of the time (specifically, his Proposition 2) proving that 
effective epimorphisms give rise to a projective class,
although he doesn't use this terminology.
It is also a bit more complicated, in the (**) case, in that a
transfinite version of Kan's $\Ex^{\infty}$ functor~\cite{ka:exi}
is required, since we make a weaker smallness assumption.
Quillen's argument in this case can be interpreted as a verification
of the hypotheses of the recognition lemma for cofibrantly generated
model categories (our Proposition~\ref{pr:rl}).

\providecommand{\bysame}{\leavevmode\hbox to3em{\hrulefill}\thinspace}

\bigskip
\bigskip

\noindent
\textbf{Notes added after publication}
\medskip

We indicate here some corrections to the published version.  We have not
updated the body of the paper except to add footnotes which point the
reader here.

\bigskip

At the end of the first paragraph of Example~\ref{ex:identity}, we
claim without proof that there is a model structure on the category of
dg modules over a dg algebra, where weak equivalences are chain
homotopy equivalences (of the underlying chain complexes) and
fibrations are degreewise split epimorphisms.
This is not at all clear, but will be explained in an upcoming
paper by Tobias Barthel, Peter May and Emily Riehl.
Thanks to Peter May and Tobias Barthel for pointing out this error.

\bigskip

In condition (**) used in Theorem~\ref{th:simp}, we should have required
that the objects of $\cS$ be $\omega$-small with respect to all
split monomorphisms in $\cA$, not just small with respect to some cardinal.
(In this form, the (**) case of Theorem~\ref{th:simp} can also be found
as Theorem~II.5.6 in Paul G.\ Goerss, John F. Jardine, Simplicial Homotopy
Theory, Birkh\"{a}user, 1999.)
This also means that the group $G$ in Corollary~\ref{co:G-sSet}
needs to be finite.
Thanks to Theo B\"{u}hler for pointing out this error.

\bigskip

Hirschhorn's book~\cite{hirschhorn} has now been published as:
Hirschhorn, Philip S. \emph{Model categories and their localizations.} 
Mathematical Surveys and Monographs, 99. 
American Mathematical Society, Providence, RI, 2003.

\bigskip
The book~\cite{dwhika} has now been published as:
Dwyer, William G.; Hirschhorn, Philip S.; Kan, Daniel M.; Smith, Jeffrey H. 
\emph{Homotopy limit functors on model categories and homotopical categories.} 
Mathematical Surveys and Monographs, 113. 
American Mathematical Society, Providence, RI, 2004.

\bigskip
\bigskip

\end{document}